\documentclass{article}
\usepackage{graphicx} % Required for inserting images
\usepackage{amsmath}
\usepackage{amssymb}
\usepackage{xcolor}
\usepackage{braket}
\usepackage{stmaryrd}
\usepackage{geometry}
\usepackage{hyperref}
\usepackage[normalem]{ulem}%

\usepackage{comment}
\hypersetup{
    colorlinks,
    linkcolor={red!50!black},
    citecolor={blue!50!black},
    urlcolor={blue!20!black}
}

\newtheorem{rem}{Remark}
\newcommand{\avg}[1]{\braket{#1}}
\newcommand\braced[1]{\left\{{#1}\right\}}

\newcommand{\pad}[2]{\frac{\partial #1}{\partial #2}}
\newcommand{\tder}[2]{\frac{{\rm d} #1}{{\rm d} #2}}
\newcommand{\GG}{{\mathcal G}}
\newcommand{\mean}[1]{\braket{#1}}
\newcommand{\order}{{\mathcal O}}
\newcommand{\fluctint}[1]{\llbracket #1 \rrbracket}
\newcommand{\xeul}{\chi}
\newcommand{\xlag}{x}
\newcommand{\xdum}{\eta}
\newcommand{\FFT}{{\rm FFT}}
\newcommand{\IFFT}{{\rm IFFT}}
\newcommand{\lep}{\eta}

\newcommand{\pmean}{\overline{p}}
\newcommand{\umean}{\overline{u}}

\newcommand{\nonlin}{{\mathcal N}}

\newcommand{\revA}[1]{{#1}}

\title{Solitary wave formation in the compressible Euler equations}
\author{David I. Ketcheson\thanks{{\texttt{david.ketcheson@kaust.edu.sa}}, Applied Mathematics and Computational Science, CEMSE Division, King Abdullah University of Science and Technology (KAUST), Thuwal, 23955-6900, Kingdom of Saudi Arabia}
\and Giovanni Russo\thanks{\texttt{giovanni.russo1@unict.it}, Department of Mathematics and Computer Science, University of Catania, Viale A.\ Doria 6, 95125 Catania, Italy}}

\begin{document}

\maketitle

\begin{abstract}
We study the behavior of perturbations in a compressible
one-dimensional inviscid gas with an ambient state consisting
of constant pressure and periodically-varying density.  We show through asymptotic analysis that long-wavelength perturbations approximately obey a system of dispersive nonlinear wave
equations.  Computational experiments demonstrate that solutions of the 1D Euler equations agree well with this dispersive model, with solutions consisting mainly of solitary waves.  Shock formation seems to be avoided for moderate-amplitude initial data, while shock formation occurs for larger initial data.  We investigate the threshold for transition between these behaviors, validating a previously-proposed criterion based on further computational experiments.  These results support the existence of large-time non-breaking solutions to the 1D compressible Euler equations, as hypothesized in previous works.
\end{abstract}

\section{Introduction}
One of the most fundamental questions about any dynamical system
is whether every solution trajectory eventually approaches a
constant state.  For nonlinear hyperbolic conservation laws in
one dimension, decay to a state that is constant in both space and time seems to be the generic behavior.
It has been known since the work of Riemann that solutions of the Euler equations of compressible gas dynamics in one dimension generically form shocks and then decay to an asymptotically constant state.  This result was made rigorous in the work of Glimm \& Lax \cite{glimm1970decay}, wherein it was shown that solutions of the Cauchy problem for a genuinely nonlinear $2 \times 2$ hyperbolic system decay in time at a rate of $1/\sqrt{t}$, while on a periodic domain they decay at a rate $1/t$.

Remarkably, for the Euler equations on a periodic domain, substantial evidence has been provided to suggest that
there exist solutions that do not form shocks or decay, and even to suggest that such solutions represent the typical long-time behavior \cite{majda1988canonical,celentano1995finite,vaynblat1996strongly,shefter1999quasiperiodic,temple2009paradigm,temple2011time,temple2023nonlinear}, \revA{or that solutions may exhibit greatly delayed and weakened shock formation \cite{hunter2019resonant}.}
These solutions seem to avoid shock formation through a resonant interaction between the nonlinear acoustic fields and the linearly degenerate entropy field.

Nevertheless, for the Cauchy problem, shock formation and entropy decay still represent the typical expected behavior.  After all, for any localized perturbation it is natural to expect that the variation in the different characteristic fields will eventually separate; once the fields are not interacting with each other, the resulting simple waves will inevitably form $N$-waves and decay.

Given these contrasting expectations for problems on the real line versus a in periodic
domain, it is interesting to consider solutions evolving on the real line but in the
presence of a periodic spatial background state.
For the 2$\times$2 hyperbolic $p$-system with spatially-periodic
coefficients (but on an infinite domain), it has been observed that some solutions exhibit
solitary waves instead of shocks, and do not decay asymptotically \cite{leveque2003,ketcheson2012shock}.
For the shallow water system with spatially-periodic
bathymetry, similar behavior has recently been reported
\cite{KLR2023}, with solutions exhibiting
solitary waves and no shock formation nor decay.
Both of these systems can formally be written 
as a larger system of conservation laws, by 
introducing an evolution equation for the spatially-varying
coefficient.  When viewed this way, each of these systems
has two genuinely nonlinear characteristic fields and
one linearly degenerate field.

The Euler equations of compressible gas dynamics have
a similar structure, with two genuinely nonlinear fields
and one linearly degenerate field.  
Indeed, for smooth solutions the Euler equations can be written in Lagrangian
coordinates as the $p$-system mentioned above: a genuinely nonlinear system of two conservation laws with a third field (the entropy) that varies in space but is constant in time \revA{(note that while the Lagrangian and Eulerian formulations of gas dynamics are 
provably equivalent \cite{wagner1987equivalence} in general,
the reduction to the $p$-system is valid only in the absence of shocks).}
It is thus natural to wonder whether 
solutions of the Euler equations might
behave in a manner similar to the aforementioned
solutions of the shallow water system or $p$-system.
Indeed, the interaction between the genuinely nonlinear acoustic fields and the linearly degenerate entropy field plays a key role in the analysis of the periodic solutions discussed above
\cite{majda1988canonical,shefter1999quasiperiodic,temple2011time,temple2023nonlinear}.

Here we conduct a detailed study of solutions of the 1D Euler equations on an unbounded domain in the presence of an initial spatially-periodic entropy variation.  We use asymptotic analysis and numerical solutions to understand their behavior.  Computational experiments suggest that the fate of these solutions depends on the size of the initial data compared to the size of the variation in the entropy.
 Unsurprisingly, if the initial data are sufficiently large then the solution can include the formation of large shocks. We thus focus on solutions whose amplitude is not too large relative to the entropy variation.  Our analysis and computational experiments provide strong evidence for the following remarkable properties of such solutions:
\begin{enumerate}
    \item There is no shock formation.
    \item Positive initial disturbances generically develop into solitary wave trains, and are accurately described by a pair of dispersive homogenized wave equations.
    \item These solitary waves persist indefinitely; thus the solution does not exhibit long-term decay.
\end{enumerate}
This behavior is just the opposite
of what is typically expected for the 1D Euler equations, but it is not unexpected based on the work of LeVeque \& Yong \cite{leveque2003}.
This does not violate
existing results on decay for the Cauchy problem, such as the one mentioned above, because those results deal only with small, localized initial
data, whereas here the initial entropy field is assumed to vary over the whole real line.
Nevertheless, our results suggest that decay of solutions should not be viewed as a universal
property of nonlinear hyperbolic systems, even on
unbounded domains.

In the remainder of this section we provide an example showing solitary wave formation and then we review the Euler equations and their connection with the $p$-system.  In Section \ref{sec:pert} we derive a dispersive model based on asymptotic analysis and the work of LeVeque \& Yong \cite{leveque2003}.  In Section \ref{sec:well-posed} we show how the dispersive model, which is ill-posed in its natural form, can be rewritten to be at least linearly well-posed.
In Section \ref{sec:comparison} we show that numerical solutions of the dispersive model agree well with solutions of the 1D Euler equations for moderate-amplitude initial data, but not for large-amplitude initial data.
In Section \ref{sec:traveling-waves} we study traveling wave solutions of the dispersive model and compare them with traveling waves obtained from numerical solution of the full 1D Euler equations.
In Section \ref{sec:entropy} we study the entropy evolution.  First we recall a criterion for shock formation proposed in \cite{ketcheson2012shock} and verify that it accurately predicts whether shocks will form in the 1D Euler setting.  We observe that for moderate-amplitude initial data the numerical entropy change becomes very small as the mesh is refined, and we show that for such data the solution can be accurately computed using a non-conservative numerical method.
Finally, in Section \ref{sec:random} we investigate the generality and applicability of these findings in a scenario where the ambient state is not perfectly periodic, but has some random modulation.
The code to reproduce most of the calculations and figures in this work is available online\footnote{\url{https://github.com/ketch/Euler_1D_homogenization_RR}}.

\subsection{A motivating example}\label{sec:motivating}
Let us consider the Cauchy problem for the 1D Euler equations (see \eqref{euler} below).
The background density is given by
\begin{subequations} \label{initial-data}
\begin{align} \label{pwc-variation}
    \hat{\rho}(\xeul) & = \begin{cases}  
        1/4 & 0 \le \xeul - \lfloor \xeul \rfloor < 1/2 \\
        7/4 & 1/2 \le \xeul - \lfloor \xeul \rfloor < 1.
        \end{cases}
\end{align}
and we define the material coordinate
$\xlag = \int \rho(\xeul)d\xeul$.
Then the initial data is given by
\begin{align} %\label{initial-data}
    p(\xlag,0) & = 1 + \frac{3}{20} \exp(-\xlag^2/25) \\
    \rho(\xlag,0) & = p(\xlag,0)^{1/\gamma} \hat{\rho}(\xeul(\xlag)) \\
    u(\xlag,0) & = 0.
\end{align}
\end{subequations}
This initial data corresponds to a gas in which the initial density and temperature vary
periodically in space, such that the pressure is constant over most of the domain.
Near $x=0$, there is an initial positive pressure perturbation that will split into
two opposite-traveling pulses.
%Note that the periodic density function
%$\hat{\rho(\xlag)}$ is chosen such that, in Eulerian coordinates, half of the domain would be
%occupied by each density of gas, if the pressure were constant.

\begin{figure}
\centering
    \includegraphics[width=6in]{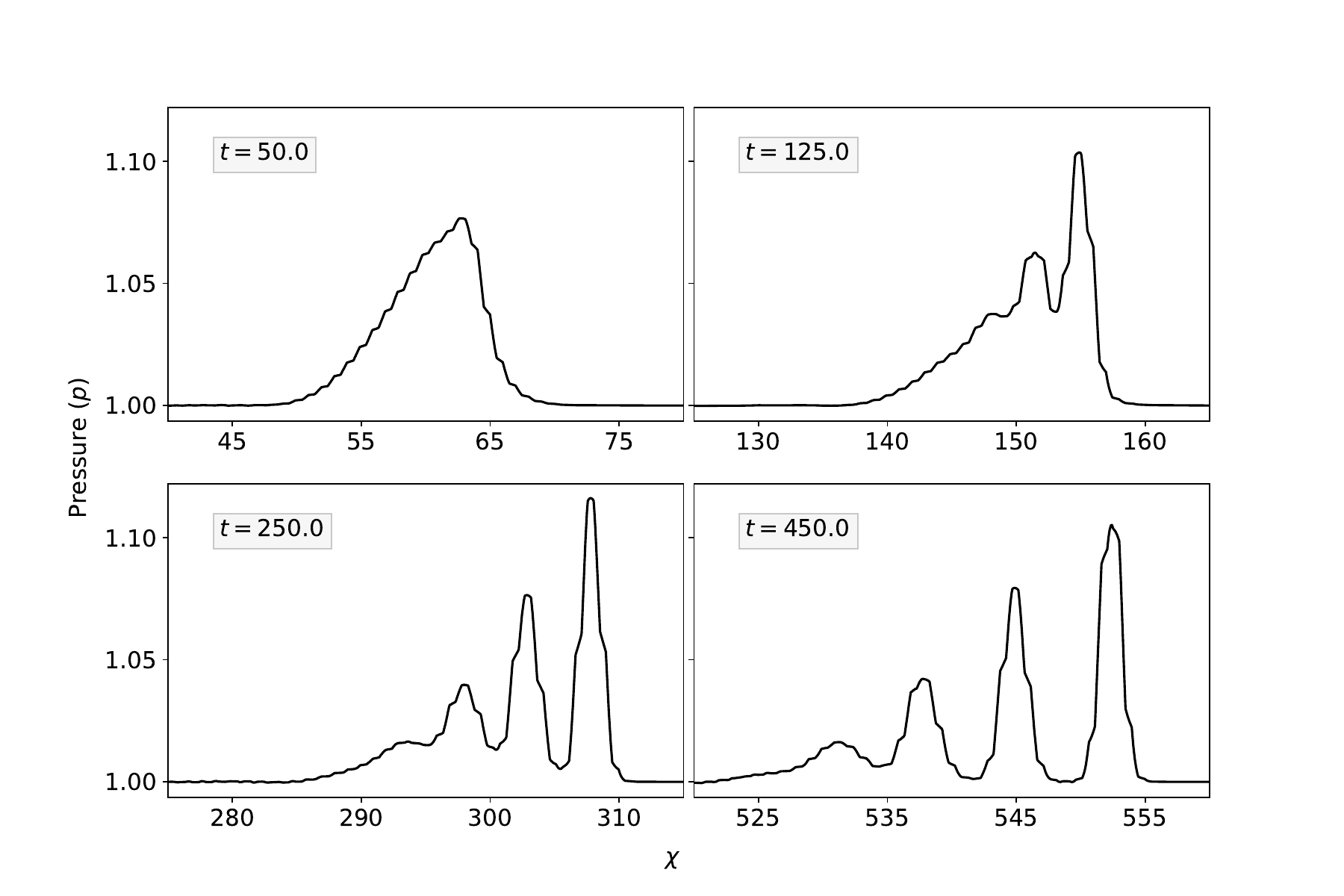}
    \caption{
 Solution of the Euler equations for a disturbance \eqref{initial-data} propagating through a periodically oscillating entropy field \eqref{pwc-variation}.}\label{fig:intro-example}
\end{figure}

Figure \ref{fig:intro-example} shows the resulting pressure solution at four different times.
The initial pulse first steepens, but at later times,
rather than the usual $N$-wave formation typical of localized or periodic solutions of 1D hyperbolic conservation laws, we observe
the emergence of a series of solitary waves, reminiscent of the behavior of dispersive nonlinear
wave equations such as Korteweg-de Vries (KdV).

\begin{figure}
\centering
    \includegraphics[width=6in]{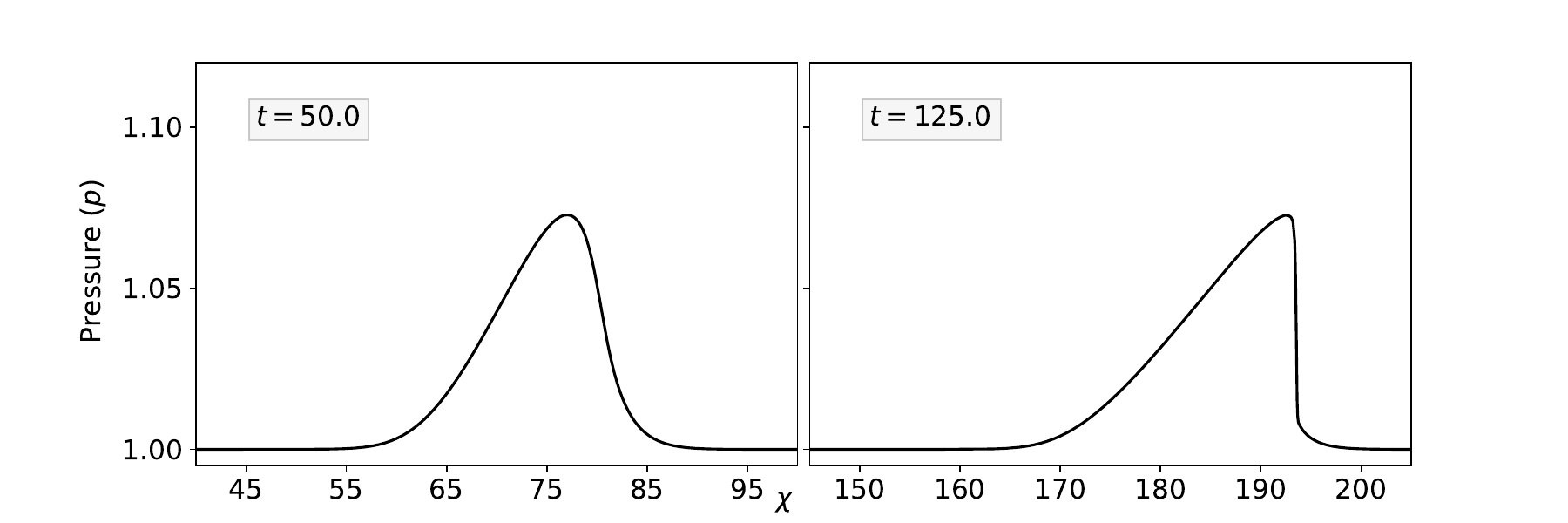}
    \caption{
 Solution of the Euler equations for a disturbance propagating through a constant entropy field.}\label{fig:intro-example-const}
\end{figure}

By way of contrast, in Figure \ref{fig:intro-example-const} we show the behavior of the same disturbance when the background density is taken to be constant and equal to the average of the oscillating density from the previous example: $\hat{\rho}(x)=1/2$.
In this case we see the $N$-wave formation, and the solution will asymptotically decay (as $t\to \infty$) to a spatially-uniform state. 
This is typical of solutions to 1D hyperbolic conservation laws.  Note also that the solution propagates significantly faster.

\subsection{Model equations}
We start from the one-dimensional Euler equations
\begin{subequations} \label{euler}
\begin{align}
    \rho_t + (\rho u)_\xeul & = 0, \nonumber \\
    (\rho u)_t +  (\rho u^2 + p)_\xeul & = 0,  \label{eq:Eulerian}\\
    \left(\frac{1}{2}\rho u^2 + \rho e\right)_t + \left( \frac{1}{2}\rho u^3 + \rho e u + u p \right)_\xeul & = 0, \nonumber
\end{align}
\end{subequations}
representing conservation of mass, momentum, and energy.
Here $\rho(\xeul,t)$ is the mass per unit volume, $u(\xeul,t)$ the gas velocity, $p(\xeul,t)$ the gas pressure, and $e$ is the internal energy per unit mass.  We consider a polytropic gas, for which
\begin{equation}\label{eq:EOS-gamma}
    e=\frac{1}{\gamma-1}\frac{p}{\rho},
\end{equation} 
where $\gamma = c_p/c_v$ is the polytropic constant, given by the ratio of specific heats respectively at constant pressure and constant volume.
The Euler equations can be conveniently written in Lagrangian form by adopting the mass coordinate:
\begin{equation}\label{eq:Euler-to-Lagrange}
    \xlag = \xlag(\xeul,t) = \xlag_0 + \int_{\xeul_0}^\xeul \rho(\tilde \xeul, t)\, d\tilde \xeul.
\end{equation}
The difference $\xlag-\xlag_0$ measures the mass of the gas (per unit area) between position $x_0$ and $x$ at time $t$.
Notice that for any $t$ the relation between the mass Lagrangian coordinate $\xlag$ and the Eulerian coordinate $\xeul$  relation is invertible. 
%The inverse relation can be expressed as 
%\begin{equation}\label{eq:Lagrange-to-Euler}
%    \xeul = \xeul_0  + \int_{\xlag_0}^\xlag v(\tilde \xlag, t)\, d\tilde \xlag
%\end{equation} 
%where $v = 1/\rho$ is the specific volume. 
%Note also that for a fixed $\xeul_0$, in general $\xlag_0$ will depend on time, and vice versa, fixing $\xeul_0$, $\xlag_0$ may depend on time. 

In the new coordinates, the Euler equations can be written as 
%In a Lagrangian frame these can be simplified to
\begin{align}
    v_t - u_\xlag & = 0, \nonumber \\
    u_t + p_\xlag & = 0, \label{eq:Lagrangian} \\
    \left(\frac{1}{2}u^2 + e\right)_t + (up)_\xlag & = 0. \nonumber
\end{align}
where $v = 1/\rho$ is the specific volume. 
%These equations are in conservation form, and provide a complete description of the compressible Euler equations 
%in one space dimension, equivalent to system \eqref{eq:Eulerian} in Eulerian coordinates.

In this paper we are interested in studying the propagation of weakly nonlinear waves.  In particular, we plan to show that such waves that propagate on an unperturbed state with constant pressure and zero velocity and periodic density will never form shocks. 
For such a reason we shall also consider a different form of the equations. 
Manipulating the energy equation, making use of the first principle of thermodynamics and of the 
first two equations, one obtains that, for smooth solutions, the entropy density $s = S/c_v$ does not depend 
on time (see, for example, \cite{whitham2011linear}, Chapter 6). The last equation of \eqref{eq:Lagrangian} can then be replaced by 
\begin{equation}
    s_t = 0 \label{eq:entropy}
\end{equation}
The entropy for a polytropic gas is given by
\[ 
    s = \log(pv^\gamma) + {\rm const},
\]
which can be written as 
\begin{equation*} \label{eq:entropy_expression}
    s - s_* = \log\left(\frac{p}{p_*}\left(\frac{v}{v_*}\right)^\gamma\right),
\end{equation*}
where $p_*$ and $v_*$ denote  a particular reference state, and  $s_*$ the  corresponding entropy density.
Solving for $p$, 
%Pressure $p$, specific volume $v$ and entropy $s$ are linked by 
the functional relation 
$p=p(\rho,s)$ becomes
\begin{equation}
    p = p_*e^{s-s_*}(v_*/v)^\gamma \label{eq:p(s)}.
\end{equation}
From Eq.~\eqref{eq:entropy} it follows that the entropy density is a function of space only, 
\[
    s = s(\xlag)
\]
In particular, $s(\xlag)$ is determined by considering an initial condition which is a perturbation of a stationary state 
\[  
    u_0 = 0, \> p_0 = p_*, \> v_0 = v_0(\xlag),
\]
so that 
\revA{\begin{equation}  
    \label{eq:s(x)}
    e^{s-s_*} = (v_0(\xlag)/v_*)^\gamma, \quad {\rm i.e.} \quad s-s_* = \gamma \log(v_0(\xlag)/v_*).
\end{equation}}

In the classical statistical theory of gases, the entropy is defined up to an additive constant. Without loss of generality we take $s_*=0$.

\begin{figure}
\centering
    \includegraphics[width=4in]{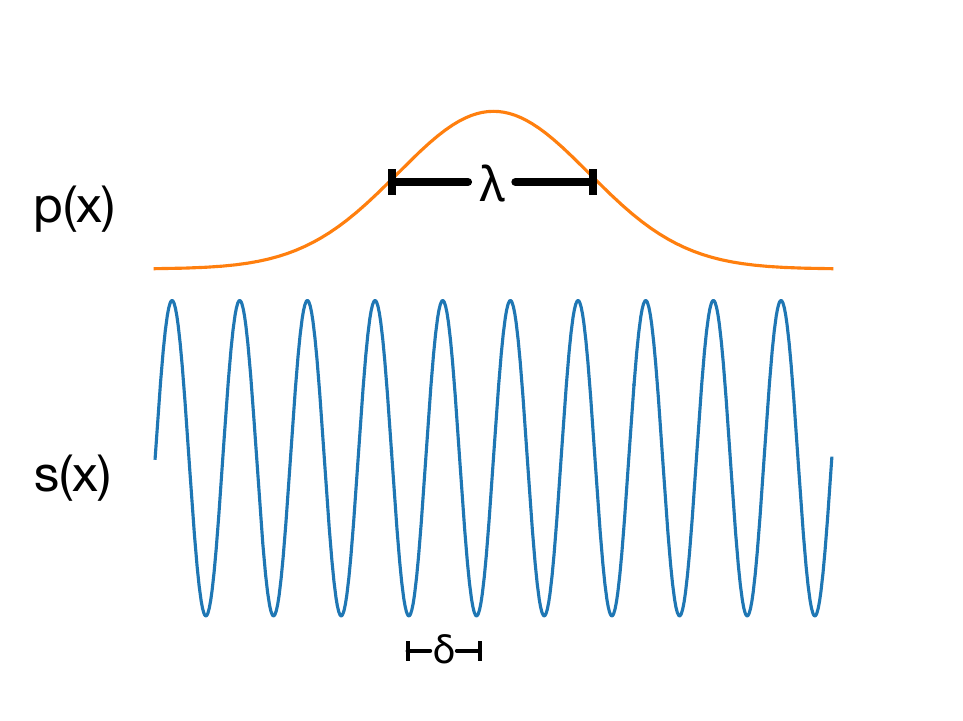}
    \caption{
 Scale of spatial variation in $S$ and $p$.}\label{fig:scaling}
\end{figure}

Using Eq.~\eqref{eq:p(s)}, and considering that $s = s(\xlag)$, the $3\times 3$ system 
\eqref{eq:Lagrangian} reduces to the $2\times 2$ system with space-dependent coefficients
\begin{subequations}
\begin{align} \label{p-system}
    v_t - u_x & = 0 \\
    u_t + p(v,s(x))_x & = 0.
\end{align}
\end{subequations}
We consider an unperturbed situation in which the initial density $v_0(\xlag)$ is a periodic function of $\xlag$, with period $\delta$\footnote{In the examples in this work we use initial data with average density equal to 1, which conveniently means that the period is the same in both the Eulerian and Lagrangian coordinates.} and assume that the initial condition is a perturbation of amplitude 
$O(\delta)$, with a typical wavelength large compared to $\delta$, see Figure \ref{fig:scaling}.

Given that the specific volume has oscillations whose amplitude is $O(1)$, we prefer to use $p$ in place of $v$ as dependent variable, since the pressure is expected to be much less oscillatory than 
$v$. 

We therefore reformulate the problem as 
\begin{align}
    \left(\pad{v}{p}\right)_{\xlag={\rm const}} p_t - u_\xlag & = 0 \\
    u_t + p_\xlag & = 0
\end{align}
Since the entropy depends only on space and not on time, we have
\[  
    \left(\pad{v}{p}\right)^{-1}_{\xlag={\rm const}}  = \left(\pad{p}{v}\right)_{\xlag={\rm const}}  
    = \left(\pad{p}{v}\right)_{s={\rm const}}  = -\frac{\gamma p}{v} = - c^2
\]
where $c$ denotes the sound speed in Lagrangian coordinates. It represents the mass (per unit surface) crossed by a small pressure perturbation per unit time. 
It is related to the classical Eulerian sound speed $c_E$ by the relation 
\[  
    c_E = v c
\]
In the field variables $(p,u) $ the system can be written as 
\begin{subequations} \label{psystem-pu}
\begin{align}
    p_t + c^2 u_\xlag & = 0, \\
    u_t + p_\xlag & = 0.
\end{align}
\end{subequations}
\revA{Making use of Eqs.~\eqref{eq:p(s)}-\eqref{eq:s(x)}, $c^2$ may be expressed as 
\[  
    c^2 = \frac{\gamma p}{v} = c_*^2e^{-s(x)/\gamma}\left(\frac{p}{p_*}\right)^{1+1/\gamma}
    = c_*^2K(\xlag)\left(\frac{p}{p_*}\right)^{1+1/\gamma} 
\]
where, for convenience,} we introduced the quantities
\begin{align}
    K(\xlag) & \equiv e^{-s(\xlag)/\gamma} = v_*/v_0(x) \\
    c_* & \equiv \frac{\gamma p_*}{v_*}.
\end{align}

In this work we consider the Cauchy problem for propagation of waves in a compressible gas in which the entropy field is initially periodic in space.  We consider the problem in Lagrangian coordinates and suppose that no shocks form, so that $s=s(x)$ and \eqref{psystem-pu} holds.
We then derive a constant-coefficient (i.e.\ effective medium) equation
that approximates the $p$-system \eqref{psystem-pu}, using multiple-scale
perturbation analysis.  An analysis of the system \eqref{psystem-pu} with a different constitutive relation and physical interpretation was
carried out by LeVeque \& Yong, who found that the resulting equations admit traveling solitary wave solutions, in good agreement with numerical experiments \cite{leveque2003}.  In that work the authors considered a special exponential equation of state and focused on piecewise constant materials.  Here we revisit the analysis in the context of compressible gas dynamics and consider general periodic variation of $s(x)$.  We also make a connection with related analysis in \cite{ketcheson2012shock} and propose a criterion for shock stability in a compressible gas with a spatially periodic entropy.

In all numerical experiments in this work, we take $v_*=p_*=1$.

\section{A Dispersive Effective Model} \label{sec:pert}

\subsection{Averaging operators}\label{sec:averaging}
In what follows we make use of certain averaging
operators, which appear in the literature but are recalled
here for the convenience of the reader.

First, the integral (or average) of $f$ over one period, denoted by $\mean{f}\in\mathbb{R}$, is defined as
\[
    \mean{f} := \int_0^1 f(y)\, dy. 
\]
\revA{For brevity we will also sometimes denote $\mean{f}$ by $\overline{f}$.}
Second, the fluctuating part of the function $f$, denoted by $\{ f \}$, is defined as 
$$\{f\}(y):=f(y)-\mean{f}.$$
Finally, the fluctuating part of the antiderivative of the fluctuating part, denoted by 
$\llbracket f \rrbracket$, is defined for any $y$ as
\revA{
\[
\fluctint{f}(y):= \left\{\int_0^y \braced{f(\xdum)}\, d\xdum\right\},
\]
that is,
\begin{equation}\label{[[]]expldef}
    \fluctint{f}(y) = \int_0^y\left\{f\right\}(\xdum)d\xdum  - \int_0^1 \int_0^\tau \left\{f\right\}(\xdum)d\xdum d\tau.
\end{equation}
}
Clearly, we have 
\[
\mean{\{f\}}=0 \quad \quad \text{and} \quad \quad \mean{\fluctint{f}}=0.
\]

Some useful properties of the $\fluctint{\cdot}$ operator, 
which will be used in this work, are provided in \cite[Appendix A]{KLR2023} .

\begin{rem}
We will often write $\braket{f}$
even for functions $f$ that are independent of $y$ \emph{a priori}, in order to emphasize which factors do not depend on $y$. 
Also, note that while $\mean{f}$ is $y$-independent, $\{f\}$ and $\fluctint{f}$ depend on $y$.
\end{rem}

\begin{rem}
\revA{Note that throughout this work, $K^{-1}$ denotes the reciprocal $1/K(x)$ and not the functional inverse.}
Furthermore, the frequently appearing average $\braket{K^{-1}}$ is equal to the total mass per period in the background state; i.e. $\mean{K^{-1}}=\rho_* \xeul(\delta)$.
\end{rem}

\subsection{Effective model}
LeVeque \& Yong \cite{leveque2003} studied the propagation of elastic 
waves in a periodically-varying medium, using the model
\begin{subequations} \label{LYmodel}
\begin{align}
    \sigma_t - K(x) G(\sigma) u_x & = 0 \\
    \rho(x) u_t - \sigma_x & = 0.
\end{align}
\end{subequations}
The $p$-system \eqref{psystem-pu} is a special case of \eqref{LYmodel}, with $\rho(x)=1$
and the notational changes
\begin{subequations} \label{correspondence}
\begin{align}
    \rho(x) & \leftrightarrow 1 \\
    u(x,t) & \leftrightarrow u(x,t) \\
    \sigma(x,t) & \leftrightarrow -p(\xlag,t) \\
    K(x) & \leftrightarrow e^{-s(\xlag)/\gamma} =: K(\xlag) \\
    G(\sigma) & \leftrightarrow  c_*^2\left(\frac{p}{p_*}\right)^{1+1/\gamma} =: G(p).
\end{align}
\end{subequations}
Here we have overloaded the definitions of $K$ and $G$ in order to facilitate
comparison between the analysis here and in \cite{leveque2003}.
The correspondence \eqref{correspondence} allows us to directly transfer results of the analysis conducted in \cite{leveque2003}
to the present setting.  To do so, we introduce a fast spatial variable $y=x/\delta$, formally independent of $x$, such that partial derivatives are transformed as
\[
    \frac{\partial}{\partial x} \to \frac{\partial}{\partial x} + \delta^{-1} \frac{\partial}{\partial y}.
\]
We suppose there exist power series for $p, u$ in terms of $\delta$:
\begin{align}
    p(x,y,t) & = p^0(x,t) + \delta p^1(x,y,t) + \delta^2 p^2(x,y,t) + \cdots \\
    u(x,y,t) & = u^0(x,t) + \delta u^1(x,y,t) + \delta^2 u^2(x,y,t) + \cdots.
\end{align}
We can also expand $G(p)$ as a power series:
\begin{align}
    G(p) & = G(p^0) + G'(p^0)(p-p^0) + \frac{1}{2}G''(p^0)(p-p^0)^2 + \cdots \\
         & = \frac{c_*^2}{p_*^{1+1/\gamma}}\left( (p^0)^{1+1/\gamma} 
          + (1+1/\gamma)(p^0)^{1/\gamma}(p-p^0)
          +\frac{1+1/\gamma}{2\gamma}(p^0)^{1/\gamma-1}(p-p^0)^2 + \cdots
         \right).
\end{align}
\revA{For brevity in the equations that follow, we now introduce the averaged quantities $\pmean(x,t) = \langle p(x,y,t)\rangle$ and $\umean(x,t) = \langle u(x,y,t) \rangle$.}
Following the analysis in \cite{leveque2003} we obtain formally
\begin{subequations} \label{homog_delta4}
\begin{align} 
\pmean_t +\frac{G}{\mean{K^{-1}}}\umean_x & + \delta^2 \frac{\mu\mean{K^{-1}}}{ G}  \left(\frac{4 \pmean_t
   \pmean_{tt} G'}{G}+\pmean_t^3 \left(\frac{G''}{G}-\frac{3
   G'^2}{G^2}\right)-\pmean_{ttt}\right)  \nonumber \\
+ \delta^4 \frac{\zeta}{\mean{K^{-1}}} &   \left(\alpha_1 \pmean_{tttt} \pmean_t - \alpha_2 \pmean_{tt} \pmean_{ttt}
   + \alpha_3 \pmean_t \pmean_{tt}^2 + \alpha_4 \pmean_t^2
   \pmean_{ttt} + \alpha_5 \pmean_t^3 \pmean_{tt}  + \alpha_6 \pmean_t^5  + \alpha_7 \pmean_{ttttt}\right) & = \order(\delta^5), \\
   \umean_t + \pmean_x & = \order(\delta^5),
\end{align}
\end{subequations} 
where $G=G(\pmean)$ and
\begin{subequations}\label{eq:list} 
\begin{align}
    \mu & = \frac{\mean{\fluctint{K^{-1}}^2}}{\mean{K^{-1}}^2} \\
    \zeta & = \mean{K^{-1}(\fluctint{K^{-1}})^2} \\
    \alpha_1 & = -9 \frac{G'}{G} \\
    \alpha_2 & =  -15\frac{G'}{G^3} \\
    \alpha_3 & = \frac{62 G'^2-\frac{37}{2} GG''}{G^4} \\
    \alpha_4 & = \frac{46 G'^2-13 GG''}{G^4} \\
    \alpha_5 & = \frac{-11 G^2 G''' - 160 G'^3 + 108 G G' G''}{G^5} \\
    \alpha_6 & = \frac{-G^3 G^{(4)}+9 G^2 G''^2+75 G'^4+13 G^2 G''' G'-\frac{165}{2} G G'^2 G''}{G^6}\\
    \alpha_7 & = \frac{1}{G^2}.
\end{align}
\end{subequations}

Note that we always have $\mu\ge 0$ and $\zeta \ge 0$.
Here we have assumed that $K(y)$ is translation-even\footnote{We say a function defined on a periodic domain is translation-even if it can be made even by a shift.}
and applied \cite[Proposition 5]{KLR2023}
in order to simplify this expression somewhat;
for a more general function $K(y)$ \eqref{homog_delta4} includes several additional terms of order $\delta^4$.
In the examples below we consider piecewise-constant or sinusoidal functions $K(y)$, both of which
are shift-periodic.
%(analogous to
%eqn. (5.17) therein)
%\begin{subequations}
%\begin{align}
%    p_t + \frac{G(p)}{\mean{K^{-1}}} u_x & = 
%     -\delta^2 \frac{\mu}{\mean{K^{-1}}} \left( G(p)u_{xxx} + 2 G'(p) p_x u_{xx} + G''(p)p_x^2 u_x \right) + \order(\delta^3) \\
%     u_t + p_x & = 0.
%\end{align}
%\end{subequations}
%where
Additionally, several terms that are present in \cite{leveque2003} instead vanish
here because of the correspondence $\rho(x)\leftrightarrow 1$.\footnote{Note also that there is a typo in \cite[eqn. (5.18)]{leveque2003}: the expressions for $C_{12}$ and $C_{22}$ are reversed.}    
\revA{As noted already in \cite{leveque2003}, these equations
possess the following symmetry: if $u(x,t)$ is a solution for $\delta=1$, then $u(x/\delta,t/\delta)$ is also a solution, for arbitrary $\delta$.  Note that under this transformation the ratio between the period of the medium and the wavelength of the solution remains the same.
In the numerical simulations below, we will simply take
$\delta=1$ and choose initial data that varies slowly relative to this scale.}
%The evolution equation for $p$ can be written, with the same order of approximation in $\delta$, as
%\[
%    p_t + \frac{G(p)}{\mean{K^{-1}}} u_x = 
%     - \delta^2 \frac{\mu}{\mean{K^{-1}}} \left( 2 G'(p) p_x u_{xx} + G''(p)p_x^2 u_x \right) - \delta^2 \mu p_{xxt} + \order(\delta^3) \\
%\]
\subsection{Alternative form of the equations}\label{sec:well-posed}
System \eqref{homog_delta4} contains high order derivatives in time. This form of the equations is not convenient for two reasons. First, the original system of PDE's is a first order hyperbolic system of two equations, therefore the initial value problem is well posed when two initial profiles are assigned. In contrast, system 
 \eqref{homog_delta4} requires four initial conditions, namely the initial values of
 $\umean$, $\pmean$, $\pmean_t$, and $\pmean_{tt}$. 

 Furthermore, system  \eqref{homog_delta4} is linearly unstable. \revA{This is a known issue with certain homogenization techniques like that employed herein; see \cite{allaire2022crime} for a general discussion, Appendix \ref{A:stability} for an analysis of the present case, and \cite{KLR2023,busaleh2024homogenized} for additional similar examples.
 What is recommended in \cite{allaire2022crime},
 and also employed in \cite{leveque2003}, is to exchange all high-order time derivatives (i.e., all except those appearing in the linear evolution term for each equation) for space derivatives, by using equality of mixed partial derivatives and keeping the same formal order of accuracy in $\delta$.
 This approach leads to a system that is also linearly unstable, but in a weaker sense -- it exhibits instability only for high-wavenumber modes, and consequently it is possible to use it for numerical simulations (as is done in \cite{leveque2003}) as long as the mesh is not too fine.}

\revA{Here we pursue an alternative approach, following what is done
in \cite{KLR2023,busaleh2024homogenized} and consistent with observations dating back to Whitham \cite[Section~13.11]{whitham2011linear}.  Namely, we
exchange some time derivatives for space derivatives by using
equality of mixed partial derivatives, but we
keep exactly one $t$-derivative in each higher-order linear term, thus obtaining a system that is linearly stable for all wavenumbers:}
\begin{subequations} \label{homog-xxt}
\begin{align}
\pmean_t + \frac{G(\pmean)}{\mean{K^{-1}}} \umean_x - \delta^2 \mu \left(\pmean_{xxt} + \frac{G'(\pmean)}{\mean{K^{-1}}} \pmean_{xx}\umean_x\right) + \delta^4 \left(\frac{\zeta}{\mean{K^{-1}}^3} - \mu^2 \right) \pmean_{xxxxt} & = \delta^4 \nonlin(\pmean,\umean) + \order(\delta^5) \label{homog-a}\\
\umean_t + \pmean_x & = 0.\label{homog-b}
\end{align}
\end{subequations}
Here $\nonlin(\pmean,\umean)$ is a function of $\pmean$ and $\umean$ and their derivatives, in which every term is nonlinear:
\begin{align*}
    \nonlin(\pmean,\umean) = & \frac{1}{\mean{K^{-1}}} \left(\frac{\zeta}{\mean{K^{-1}}^3} - \mu^2 \right)
    \Bigg[ 
        \beta \left(2\frac{G'}{G}\pmean_{xx} - \frac{G''}{G} (\pmean_x)^2\right) - 6G' \pmean_{xxx}\umean_{xx} - G'\pmean_{xx}\umean_{xxx} - 6\frac{(G')^2}{G}(\pmean_x)^2\umean_{xxx} \\
        & - \left(8  \frac{(G')^2}{G} + 6 G''\right)\pmean_x\pmean_{xx}\umean_{xx}  - 6 G'' \pmean_x \pmean_{xxx} \umean_x - 6 \frac{G' G''}{G}(\pmean_x)^3 \umean_{xx} \\
        & + \frac{\beta}{\mean{K^{-1}}} \left( \frac{G''}{G}-2(G')^2\right) (\umean_x)^2 + \frac{1}{\mean{K^{-1}}}\left(2(G')^2-6 G G''\right)\umean_x (\umean_{xx})^2 \\
        & + \frac{1}{\mean{K^{-1}}}\left(2(G')^2 - G G''\right)(\umean_x)^2\umean_{xxx} +  \frac{1}{\mean{K^{-1}}}\left(2\frac{(G')^3}{G} - G' G''\right) \pmean_{xx} (\umean_x)^3 \\
        & + \left(-9\frac{G'G''}{G} - 3 G^{(3)}\right)(\pmean_x)^2\pmean_{xx}\umean_x - 2\frac{G' G^{(3)}}{G}(\pmean_x)^4 \umean_x \\
        & + \frac{1}{\mean{K^{-1}}}\left(4\frac{(G')^3}{G} - 4 G' G'' - 6 G G^{(3)}\right) \pmean_x (\umean_x)^2 \umean_{xx} \\
        & + \frac{1}{\mean{K^{-1}}}\left(2\frac{(G')^2 G''}{G} - 2(G'')^2 - 2 G' G^{(3)} -G G^{(4)} \right) (\pmean_x)^2 (\umean_x)^3 + G' \pmean_{xxxx}\umean_x \\
        & - 2 G' \pmean_x\umean_{xxxx} - 2 \frac{(G')^2}{G}(\pmean_{xx})^2\umean_x
    \Bigg]
\end{align*}
%where for shortness we set
%\begin{align*}
%    + \delta^4 \frac{\zeta G'' }{2\mean{K^{-1}}^4}(\pmean_{xx})^2\umean_x 
%\beta & = G\umean_{xxx} + G' \pmean_{xx}\umean_x + 2G'\pmean_x\umean_{xx} + G'' 
% \umean_x(\pmean_x)^2.
%\end{align*}
The number of terms appearing here in the evolution equation for $p$ is much larger than in \cite{leveque2003}, because in that work the equation of state was chosen such that $G''=0$.
On the other hand, the evolution equation for $u$ here is much simpler, due to the fact that the variable coefficient $\rho(x)$ present in \cite{leveque2003} is replaced by a constant here.

\subsection{Linear stability}\label{sec:linear_stability}
In this section we study the well-posedness, for small initial data, of the initial value problem for the model 
equations derived in the previous section.  As mentioned above, the linearization of system \eqref{homog_delta4} is unstable for all wave numbers, and the initial value problem is ill-posed (see Appendix \ref{A:stability}).

Next we consider the system \eqref{homog-xxt}.
We linearize around an equilibrium configuration $(\umean,\pmean) = (0,p_*)$. 
Neglecting terms of $O(\delta^5)$, the linearized equations read
\begin{subequations}\label{eq:linearized}
\begin{align}
    \umean_t + \pmean_x & = 0\\
    \pmean_t + c^2 \umean_x - \delta^2\mu \pmean_{xxt} + \delta^4\nu\pmean_{xxxxt} & = 0
\end{align}
\end{subequations}
where 
\begin{equation}\label{eq:nu}
    c^2 \equiv \frac{G(p_*)}{\mean{K^{-1}}}, \quad
    \nu \equiv \frac{\zeta}{\mean{K^{-1}}^3}-\mu^2
\end{equation}
We look for solutions of the system \eqref{eq:linearized} of the form 
\[  
    \umean = \hat{u}e^{i(kx-\omega t)}, \quad \pmean = e^{i(kx-\omega t)}
\]
Plugging this {\em ansatz\/} into system \eqref{eq:linearized} we obtain 
a homogeneous system for $\hat{p}$ and $\hat{u}$:
\begin{align*}
    -i\omega\hat u + i k \hat p & = 0\\
    - i \omega \hat p + i c^2 k \hat u 
    - i \delta^2 k^2\omega \mu \hat p 
    - i\delta^4 \nu k^4 \omega \hat p& = 0.
\end{align*}
Non-trivial solutions of the above system exist provided the determinant of the coefficient matrix of the system vanishes. This gives the following
dispersion relation
\begin{equation}\label{eq:dispersion}
    c^2 k^2 - \omega^2 (1+\mu\delta^2 k^2 + \nu \delta^4 k^4) = 0,
\end{equation}
which can be solved for $\omega$:
\[
    \omega = \pm \frac{c k}{\sqrt{1+\mu\delta^2 k^2 + \nu \delta^4 k^4}}
    = \pm ck\left(1-\frac{1}{2}\mu \delta^2 k^2 + \frac{1}{8}(3\mu^2 - 4\nu)\delta^4 k^4 + \order(\delta^6 k^6)\right).
\]
For all the profiles $K(x)$ we have tested, $\nu>0$, and we conjecture that $\nu>0$
in general, which means that the systems admits linearly dispersive waves for all wave numbers $k\in\mathbb{R}$.

\subsection{Comparison of the models}\label{sec:comparison}

We now compare solutions obtained by solving the Euler equations, the 
$p$-system, and the homogenized equations.  The Euler equations in Eulerian form and the
variable-coefficient $p$-system are solved using Clawpack \cite{mandli2016clawpack}, specifically with the high-order
SharpClaw algorithm based on fifth-order WENO reconstruction and 
a ten-stage, fourth-order strong stability preserving Runge-Kutta method 
\cite{ketcheson2008highly,KetParLev13,pyclaw-sisc}.  The homogenized equations
\eqref{homog-xxt}, including terms up to $\order(\delta^4)$ are solved using a Fourier pseudospectral collocation method
in space and the three-stage, third-order strong stability preserving Runge-Kutta method of Shu \& Osher \cite{shu1988efficient}.

\subsubsection{Moderate-amplitude initial data}\label{sec:moderate}
We first consider exactly the scenario from Section \ref{sec:motivating}, with
the initial data given by 
\eqref{initial-data}.
For the pseudospectral simulation we use the domain $[-600,600]$
with periodic boundary conditions, but the simulation ends before
there is significant interaction of the waves with the boundary.  For the finite volume simulations we use the domain $[0,600]$ and
impose a reflecting (solid wall) boundary condition at $x=0$.

\begin{figure}
\centering
    \includegraphics[width=7in]{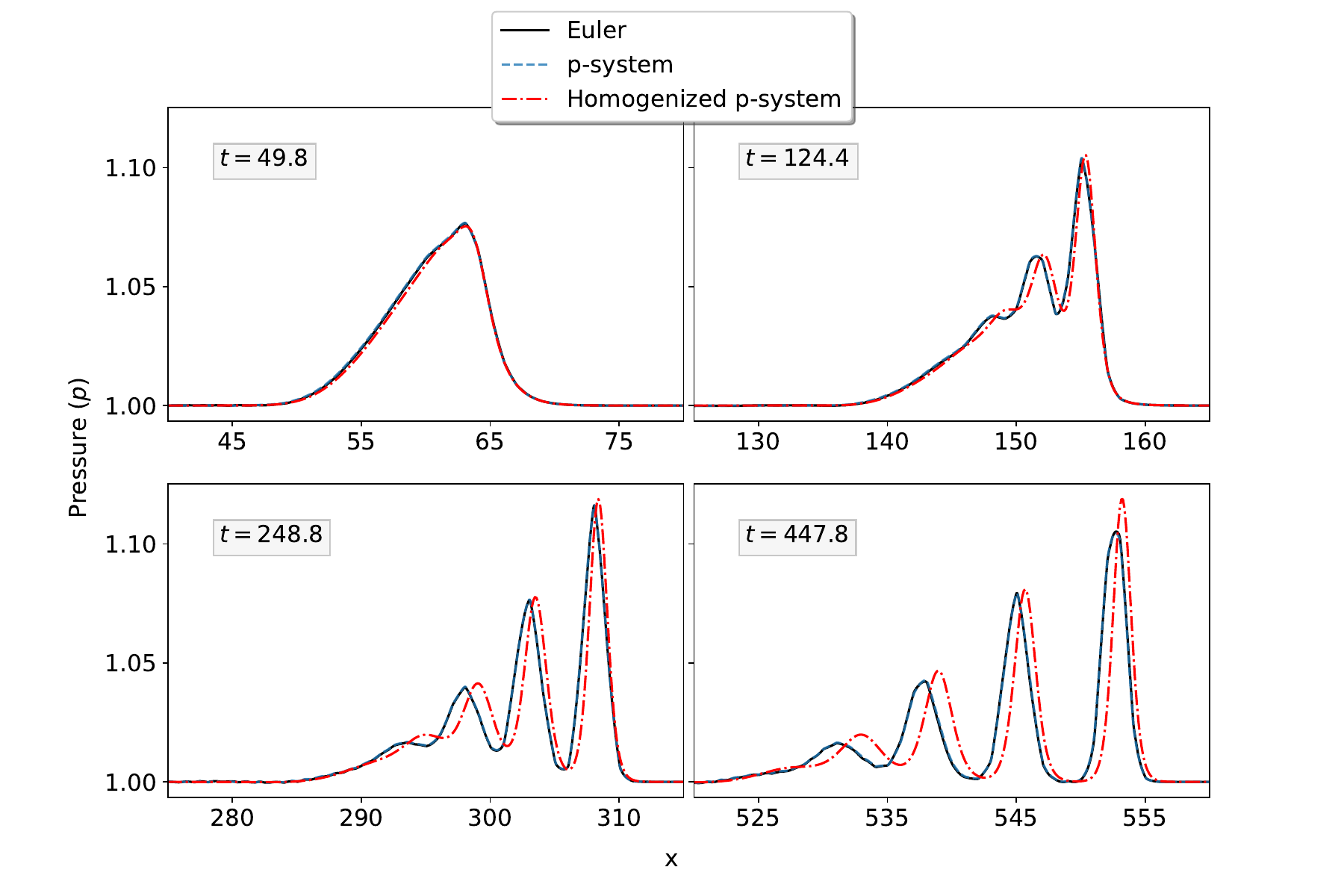}
    \caption{
 Comparison between the direct solution of the Euler equations \eqref{euler}, the $p$-system \eqref{psystem-pu}, and the homogenized equations \eqref{homog-xxt}.  The setup is the same as for Figure \ref{fig:intro-example}, but note that here the x-axis is the material coordinate so the solutions look smoother.
 The Euler and $p$-system solutions are indistinguishable, consistent with the hypothesis that shocks do not form.}\label{fig:compare3}
\end{figure}

Results are shown in Figure \ref{fig:compare3}.
The Euler and $p$-system solutions are visually indistinguishable
and seem to converge to the same values, as would be expected if shocks do not form (so that there is no temporal change in the entropy field).
The homogenized solution shows close agreement at early times, with
increasingly noticeable differences at later times.
\revA{The agreement of the homogenized solution with the detailed numerical simulation is quite remarkable, considering that the width of the pulses is approximately five times the period of the initial perturbation in the density.}

\subsubsection{Large-amplitude initial data}
\begin{figure}
\centering
    \includegraphics[width=6in]{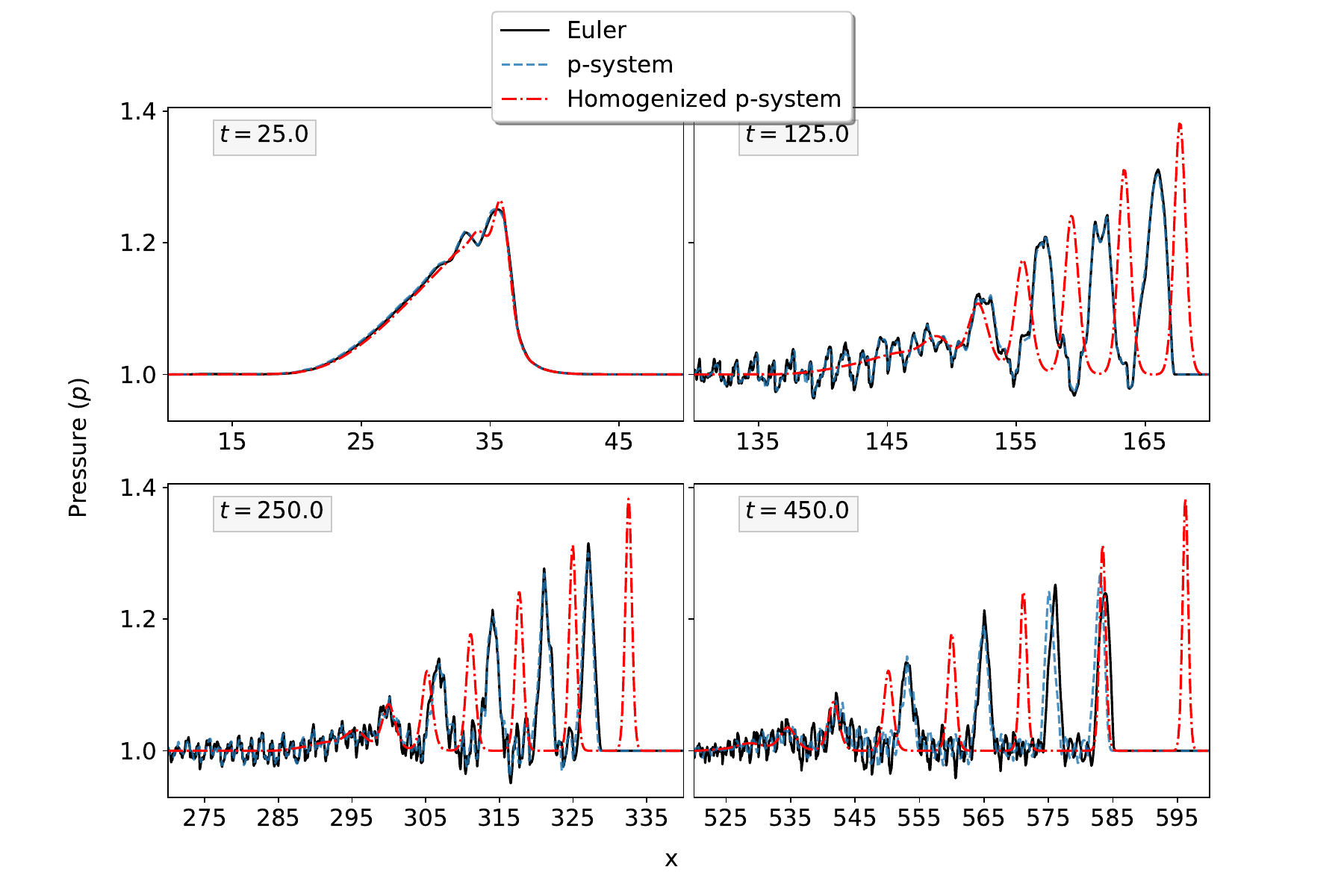}
    \caption{
 Comparison between the direct solution of the Euler equations \eqref{euler}, the $p$-system \eqref{psystem-pu}, and the homogenized equations \eqref{homog-xxt}.  The setup is the same as for Figure \ref{fig:compare3}, but with an initial pulse more than 3 times as large.
 Due to shock formation, the homogenized solutions differs markedly from the finite volume solutions.  At late times, the Euler and $p$-system solutions are visibly different. }\label{fig:compare3_shock}
\end{figure}

The behavior seen in the last example is typical for initial data that is not too large relative to the variation in the background entropy.  For larger initial data (or smaller entropy variation), typical solutions involve shock formation as is usually expected in solutions of the Euler equations.
To demonstrate this, we repeat the experiment of the previous section with one change: the amplitude of the initial pressure pulse is increased from $3/20$ to $1/2$:
\begin{align}
    p(\xlag,0) & = 1 + \frac{1}{2} \exp(-\xlag^2/25).
\end{align}
Results are shown in Figure \ref{fig:compare3_shock}.  We see that the Euler solution differs markedly from the homogenized solution and exhibits high-frequency oscillations, presumably resulting from the interaction of shock waves and the spatially-varying density.
The $p$-system solution remains close to the Euler solution with visible differences only after a very long time. 
It is well known that for weak shocks the entropy produced by the shock is very small. Indeed it is 
of $\order(([p]/p_*)^3)$, where $[p]$ denotes the jump in the pressure and $p_*$ 
is the pressure in the unperturbed region ahead of the shock 
(see \cite{whitham2011linear}, Chapter 6). 
Isentropic approximation works well even for shocks of moderate strength.
This has been observed in several contexts, 
including multilayered fluids (see \cite{phan2023numerical}).

\subsection{Exploiting the regularity of the solution}\label{sec:pseudospectral}
In the previous section we showed numerical evidence that under certain circumstances, nonlinear waves in the Euler equations do not break into shocks, and no entropy is produced. This suggests that the solution remains smooth for all times, presumably maintaining the regularity of the initial data. 
If that is the case, accurate numerical solutions do not require the use of a shock capturing scheme, and discretizations can be based on a non-conservative form of the equations.
In the case of smooth initial data it is possible to adopt, for example, a pseudospectral discretization, similar to the one adopted for the numerical solution of the homogenized equations. 

Here we write the Euler system as 
\begin{subequations}\label{eq:Euler_Eul_form}
\begin{align}
    \rho_t & = - (\rho u)_\xeul \\
    u_t & = -u u_\xeul - p_\xeul/\rho, \\
    p_t & = -u p_\xeul - \gamma p u_\xeul.
\end{align}
\end{subequations}
All spatial derivatives are computed by Fourier pseudospectral differentiation; for example, $u_\xlag \approx \IFFT(i k \FFT (u))$.
The system is integrated in time using the standard four-stage, fourth-order Runge-Kutta method.
In the computation we use a Courant number 0.9. 
As usual, we integrate the equations in the domain $[-L,L]$ with periodic boundary conditions, and plot the solution only in the interval $[0,L]$. Note that if the initial density and pressure are even functions, and the initial velocity is an odd function, then the symmetry of the solution is maintained for all times. 

As an example we compare the solution obtained by the high-order finite volume method described in the previous section, with the one obtained by the pseudospectral method.  We take as initial condition
\begin{subequations} \label{smooth-IC}
\begin{align}
    \rho(\xeul,0) & = 1 + 0.8 \cos(2\pi \xeul) \\
    u(\xeul,0) & = 0 \\
    p(\xeul,0) & = \frac{3}{20}\exp(-\xeul^2/16).
\end{align}
\end{subequations}

\begin{figure}
\centering
    \includegraphics[width=\textwidth]{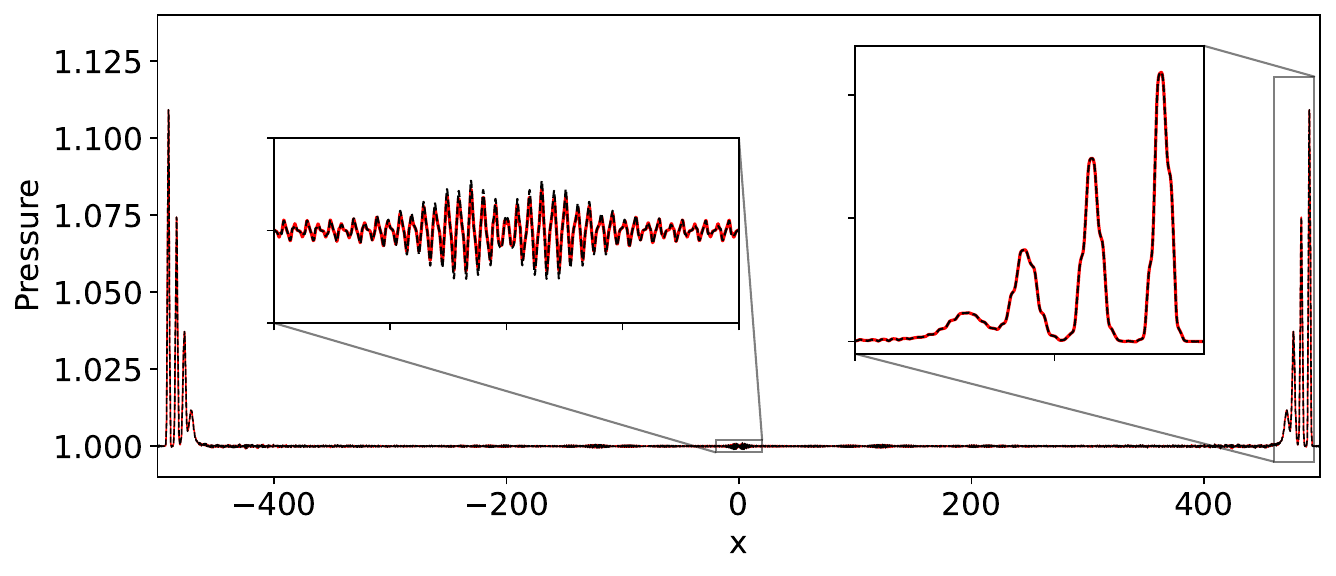}
    \caption{Pressure profile at final time obtained by Clawpack (solid red line) and Fourier pseudospectral method in primitive variables (black dashed line). Zoom insets are included to show the close agreement even for fine details of the solution.}
    \label{fig:compare_CLAW_FS}
\end{figure}

The result of the comparison is reported in Fig.~\ref{fig:compare_CLAW_FS}. 
In order to get a fully resolved calculation we used 144,000 cells with Clawpack in the interval $[-600,600]$ and only 32K points with the pseudospectral code in the domain [-512,512].
It is not our intent here to make a detailed performance comparison, and the solvers are implemented in different languages, but we remark that the pseudospectral simulation takes about one tenth of the time of the FV simulation (about 10 minutes versus about 100 minutes), with both running on a M1 Max Macbook Pro with 32 GB of RAM.

\section{Approximate traveling waves}\label{sec:traveling-waves}
In this section we compute approximate traveling wave solutions of the homogenized system \eqref{homog-xxt}.
We look for solutions which depend only on the variable 
$\xi = x-Vt$, where $V$ is the traveling speed of the wave; i.e.\ we look for solutions of the form 
$p = p(\xi)$ and $u= u(\xi)$. To simplify the notation we write simply $p, u$ in place of $\bar{p},\bar{u}$.

%\subsection{Traveling waves to $\order(\delta^2)$}
We first consider traveling waves for the approximate system in which we neglect  terms of $\order(\delta^4)$.
Inserting the traveling wave {\em ansatz} for $u$ and $p$ in system 
\eqref{homog-xxt}, 
we obtain the system of ODEs
\begin{subequations} \label{travel_system}
\begin{align}
-V p' + \frac{G(p)}{\mean{K^{-1}}} u' - \delta^2 \mu \left(-V p''' 
+ \frac{G'(p)}{\mean{K^{-1}}} p''u'\right) & = 0\\
-Vu' + p' & = 0.
\end{align}
\end{subequations}
From the second equation we deduce that $p'=V u'$, and therefore $p = p_* + V u$, since we consider propagation of traveling perturbations of the state $p=p_*$ and $u=0$.

\begin{figure}
\centering
    \includegraphics[width=0.9\textwidth]{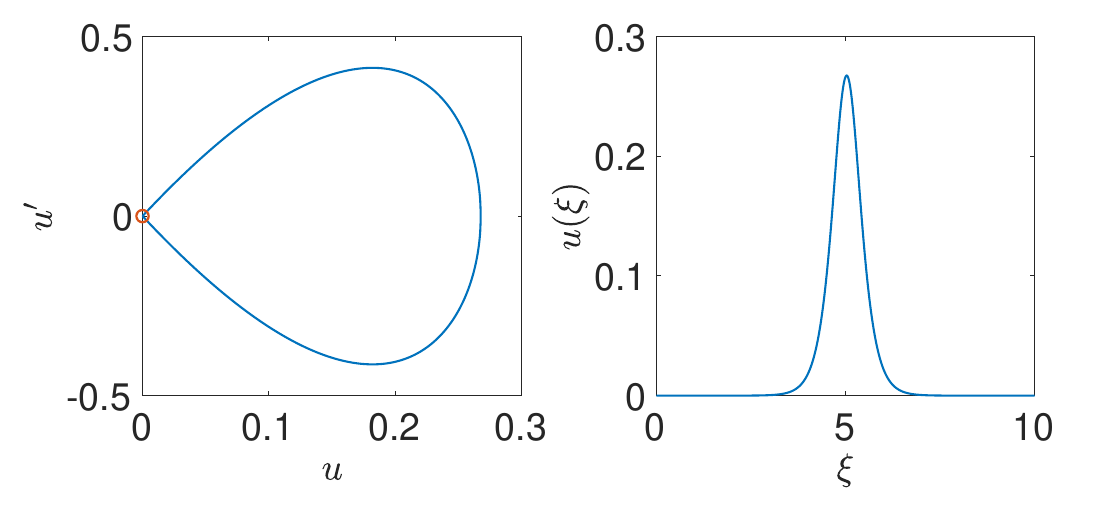}
    \caption{
 Computation of a typical traveling wave solving system \eqref{eq:uv}. Here we chose $V=1.2$. Left panel: separatrix in phase space. Right panel: traveling wave.} 
 \label{fig:travel}
\end{figure}

Inserting the dependence of $p$ on $u$ in the first equation, we obtain a third order ODE for $u(\xi)$. 
In order to integrate this equation further, we
approximate $G'(p)\approx G'(p_*)$, yielding
\begin{equation}\label{eq:u3}
-V^2 u' + \frac{G(p(u))}{\mean{K^{-1}}} u' - \delta^2 \mu \left(-V^2 u''' 
+ \frac{G'(p_*)}{\mean{K^{-1}}} V u''u'\right)  = 0.
\end{equation}
Now let $\GG(p)$ be a primitive of $G(p)$. Then ${\rm d}\GG/{\rm d}\xi = G(p(u))p' = G(p(u))Vu'$.
Using this relation, \eqref{eq:u3} can be written as 
\[
    \tder{}{\xi}\left[ -V^2 u + \frac{\GG(p(u))}{V \mean{K^{-1}}} - V^2 u + \delta^2 \mu V^2 u'' - \delta^2 \mu \frac{G'(p_*)}{\mean{K^{-1}}}V\frac12(u')^2\right] = 0
\]
which indicates that the quantity in square brackets is constant. Given that at infinity $u(\xi)$ vanishes with its derivatives, we obtain a second-order equation for $u$:
\begin{equation}
    u'' = \frac{G'(p_*)}{2V\mean{K^{-1}}}(u')^2
    -\frac{\GG(p(u))-\GG(p_*)}{\delta^2\mu V^3 \mean{K^{-1}}} 
    + \frac{u}{\delta^2\mu}.
\end{equation}
The second-order equation can be written as a first-order system of the form
\begin{equation}\label{eq:uv}
    u' = v, \quad v' = F(u,v).
\end{equation}
It turns out that $(0,0)$ is an equilibrium point for system \eqref{eq:uv}. 
The linearization around the origin reads
\[
    u' = v, \quad v' = \beta u
\]
with 
\[  
    \beta = \left(1-\frac{G(p_*)}{V^2 \mean{K^{-1}}}\right)\left(\delta^2\mu\right)^{-1}
\]
When $\beta>0$, there are two real roots, $\lambda_{\pm} = \pm \sqrt{\beta}$, and the origin is a saddle point. Integrating system 
$\eqref{eq:uv}$ with an initial condition very close to the origin, aligned with the eigenvector corresponding to the positive eigenvalue, one obtains a good approximation of a traveling wave. 
Figure \ref{fig:travel} represents a typical traveling wave.

We now use this solution as the initial guess for Newton's method in order to find a more accurate traveling wave solution.  Given an initial guess $u^{(0)}$ we iteratively solve
\begin{align*}\label{eq:SNewton}
    J[u^{(k)}] \delta u^{(k)} & = - F[u^{(k)}] \\
    u^{(k+1)} & = u^{(k)} + \delta u^{(k)}.
\end{align*}
where $F$ is the residual of the traveling wave equation:
\[
    F[u] = -u''+\alpha_0u+\alpha_1 \frac12(u')^2-\alpha_2(\GG(p_*+Vu)-\GG(p*))
\]
with 
$\alpha_0 = (\delta^2\mu)^{-1}$, 
$\alpha_1 = \frac{G'(p_*)}{V\mean{K^{-1}}}$, 
$\alpha_2 = (\delta^2\mu V^3\mean{K^{-1}})^{-1}$,
and $J$ is the Jacobian of F:
\begin{align*}
    J[u]\delta u & =  \frac{\delta F}{\delta u} \delta u \\ 
                 & = (\alpha_0-\partial^2_\xi)\delta u - \alpha_2 G(p_* + V u) V \delta u + \alpha_1 u'\delta u' \\
                 & = (\alpha_0-\alpha_2V G(p_*+Vu) - \partial^2_\xi)\delta u - \alpha_1 u''\delta u \\
                 & = (\alpha_0-\alpha_2V G(p_*+Vu) - \alpha_1 u'' - \partial^2_\xi)\delta u =: J[u]\delta u.
\end{align*}
Here we used the relation $u'\delta u' = - u''\delta u$ given that we assume that $u$ and its derivatives vanish at infinity.
We use a three-point centered difference to approximate $u''$.
%Space discretization, needed to compute $u''$ and to discretize the second derivative, can be obtained, for example, by second-order finite difference discretization. 
The Newton iteration is terminated when the residual $F[u]$ falls below a prescribed tolerance.
%The procedure is repeated until the estimate of the error becomes smaller than a desired tolerance. 

\subsection{Traveling waves to $\order{(\delta^4)}$}
We now look for traveling waves for system  
\eqref{homog-xxt}, keeping the term proportional to $\pmean_{xxxxt}$, and neglecting only the right hand side of 
\eqref{homog-a}. Let us denote by $\nu$ the coefficient of
$\pmean_{xxxxt}$, as in \eqref{eq:nu}.
Using as before the approximation $G'(p)\approx G'(p_*)$ in 
system \eqref{homog-xxt}, 
and proceeding as above, we obtain the following fourth-order equation for $u$: 
\begin{equation}\label{eq:fourth}
F[u] :=   \frac{\delta^2\nu}{\mu} u^{''''} - \left(u'' - \frac{G'(p_*)}{2V\mean{K^{-1}}}(u')^2
    +\frac{\GG(p(u))-\GG(p_*)}{\delta^2\mu V^3 \mean{K^{-1}}} 
    - \frac{u}{\delta^2\mu} \right) = 0.
\end{equation}

%The equation could be written as a first-order system with 4 equations, however it would be very difficult to capture the solution corresponding to a traveling wave by solving a suitable initial value problem. 
%For this reason, we adopt a different strategy to compute the traveling wave. We write the equation as 
%\[
%    F[u] = 0
%\]
%where $F$ may contain also differential operators acting on $u$.
%We describe the technique for the second-order model. The same can be applied to compute the traveling wave for the fourth-order model. 

We then apply the Newton method described above to \eqref{eq:fourth}, using the traveling wave computed from the $\order(\delta^2)$ approximation as initial guess.
%We start from the traveling wave computed by integrating system \eqref{eq:uv}, as describe above, and call it $u^{(0)}$. 

%A similar procedure can be adopted to compute the traveling wave corresponding to system 
%\eqref{homog-xxt}, neglecting the right hand side. 
%As initial guess for the $\order(\delta^4)$ model 
%we adopt the traveling wave of the $\order(\delta^2)$ model.

\begin{figure}
\centering
    \includegraphics[width=0.6\textwidth]{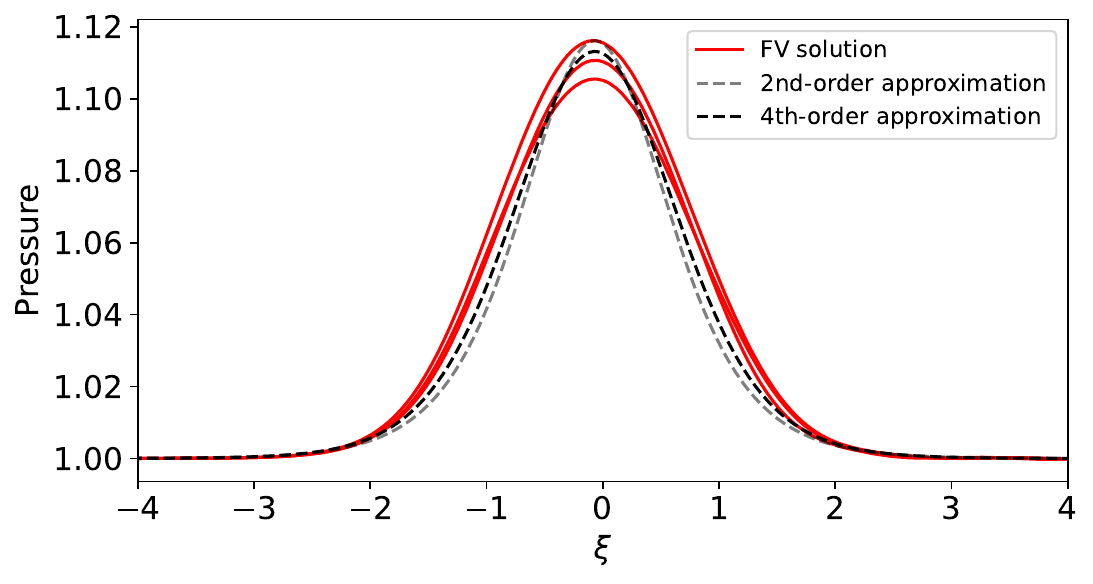}
    \caption{
 Traveling waves computed by approximately solving the homogenized equations with a traveling-wave {\em ansatz}, compared to measurements of the traveling waves observed in the full PDE simulation.  The red lines are measurements showing the variation in time at three different fixed points in space.  The dashed lines are the 2nd-order and 4th-order homogenized approximations.}
 \label{fig:compare_TW}
\end{figure}

%The amplitude of the traveling wave strongly depends on the parameter $V$, i.e. the  speed of the traveling wave. For speeds just above $c_*$ (say for $V=1.2$, while in our case $c_* \approx 1.183216$), as in Figure \ref{fig:travel}, the two approximations give almost the same traveling wave.
Figure \ref{fig:compare_TW} shows the shape of the traveling waves corresponding to the second and to the fourth-order model, respectively, with traveling velocity $V=1.222$.  We also include for comparison a traveling wave with this velocity computed by solving the initial value PDE with Clawpack.  For the latter, we show the solution as a function of $\xi$ at fixed points in space, for three different points (since the shape of the wave varies depending on the spatial location).

%\begin{figure}
%\centering
%    \includegraphics[width=0.6\textwidth]{Figures/solitons_2_4.pdf}
%    \caption{
% Traveling waves corresponding to the second and fourth-order model. Here we chose $V=1.4$.} 
% \label{fig:compare_TW}
%\end{figure}

\begin{rem}
Sometimes it is also possible to compute approximate traveling waves using Petviashvili's fixed-point iteration.
%An interesting technique for the computation of the traveling wave of various dispersive equations has been developed by Petviashvili, and is described, for example in \cite{alvarez2014petviashvili}. 
This approach is based on separating the ODE that defines the traveling wave into a linear part, which includes the term with the highest-order derivative, and a non-linear part:
\[
    L u = N(u).
\]
%After space discretization obtained, for example, by finite difference, an ingenious iterative technique is then developed, which converges to the solution under suitable assumptions. 
One of the assumptions of the procedure is that the non-linear operator is homogeneous of a certain degree. Unfortunately, this is not the case for our problem. 
%However the separation is very useful to estimate the discretization error of a given approximate solution. Indeed, if we denote by $(u,p)$ the solution of system \eqref{homog-xxt} and by $\tilde{u},\tilde{v}$ the solution of system \eqref{eq:uv}, we can compute the discretization error as $L \tilde{u}-N(\tilde{u})$. The relative
%$L^\infty$ discretization error is then 
%\[  
%    {\rm de}(\tilde{u}) = \frac{\| L \tilde{u} - N(\tilde{u})\|_\infty}{\|L\tilde{u}\|_\infty}
%\]
%Using a space discretization in $\xi$ which is fine enough that the error is mainly due to the approximation introduced in
%Eq.~\eqref{eq:u3} (i.e.\ replacing $G'(p)$ with $G'(p_*)$ as in Eq.~\eqref{eq:u3}) we obtain for the example illustrated in Figure a relative $L^\infty$ discretization error smaller than 
%$2\times 10^{-5}$ for the $V=1.2$ traveling wave.

Nevertheless, by observing  that $G(p) = p G'(p)/(1+1/\gamma)$, and approximating
$G(p)$ by $p G'(p_*)/(1+1/\gamma)$ in system \eqref{travel_system}, one can rewrite the traveling wave equation as $Lu = N(u)$, with $N(u)$ homogeneous operator of degree 2, therefore allowing the use of the Petviashvili method for the construction of the approximate traveling wave. The results are very similar to the ones obtained by Newton's method.
This was pointed out to the authors by Giuseppe Virgilio Minissale.
\end{rem}

\section{Shock formation and entropy evolution}\label{sec:entropy}
\revA{Visual inspection of the solutions computed in Section \ref{sec:moderate} suggests the possibility that no shocks are formed.  The accuracy of this statement cannot be determined from
numerical experiments alone, but in this section we attempt to
investigate it.}

A more precise method for detecting shock formation in the numerical solution is to study
the total entropy change over time.  In the exact solution of the Euler equations, the total
entropy is constant in the absence of shock formation.  For numerical solutions, we expect
some entropy change due to numerical errors.  We compute a measure of the total entropy at a 
given time from the discrete solution via
\begin{align}\label{eq:total_entropy}
    \int s(x,t_n) dx =
    %\approx \Delta x \sum_j \log(p_j^n/(\rho_j^n)^\gamma). \\
    \int s(\xeul,t_n) \rho(\xeul) d\xeul \approx \Delta \xeul \sum_j \rho_j^n \log(p_j^n/(\rho_j^n)^\gamma).
\end{align}
In order to reduce the entropy change caused by the numerical discretization, for this
test we again use the smooth initial data \eqref{smooth-IC}.
We take the domain $x\in[-256,256]$ and impose periodic boundary conditions to avoid any entropy change due to boundary effects.
We run to final time $t=200$, which
is many times greater
than the time of shock formation for a similar initial condition in a constant initial
entropy field.  In Table \ref{tab:entropy-change} we report the change in entropy:
\begin{align} \label{entropy-change}
    \Delta s_\textup{total} := %\int s(x,200) dx - \int s(x,0) dx. \\
    %\Delta S_\textup{total} := 
    \int s(\xeul,200) \rho(\xeul) d\xeul - \int s(\xeul,0) \rho(\xeul) d\xeul.
\end{align}
%\begin{align} \label{entropy-change}
%    \frac{\int S(\xeul,200) d\xeul - \int S(\xeul,0) d\xeul}{|\int S(\xeul,0) dx|}.
%\end{align}
This value should be zero in the absence of shock formation, or positive in the presence of shocks.

For this comparison we computed solutions with the high-order SharpClaw (SC) method in Clawpack
%, the second-order Lax-Wendroff-LeVeque
%algorithm originally implemented in Clawpack \cite{leveque2002finite}, 
and with the pseudospectral method
described in Section \ref{sec:pseudospectral}.
There is a small amount of entropy production in the Clawpack solutions, and in the pseudospectral solution on a coarse grid.  This entropy production decreases at approximately the expected rate of convergence for Clawpack, while the pseudospectral method shows a very small entropy production on finer grids.
This suggests that if any shock formation occurs, the shock(s) must be very weak.

%\red{When computing the evolution with a random background the relative change in the entropy at time $t=180$ has been $4.39\times 10^{-9}$.}

%\begin{table}
%\centering
%\begin{tabular}{l|rrr}
%$\Delta \xeul$ & LWL & Clawpack & Pseudospectral \\ \hline
%1/16  &         & -4.79e-3 & -2.76e-7 \\
%%1/24  &         & -7.51e-4 & -5.88e-10 \\
%1/32  &         & -1.95e-4 & 2.62e-10 \\
%%1/48  &         & -2.80e-5 & 2.02e-10 \\
%1/50  & 1.92e-2 & -2.31e-5 & 1.92e-10 \\
%%1/64 & & & 1.49e-10 &              \\
%1/100 & 4.74e-3 & -9.23e-7 & 1.10e-10 \\
%1/200 & 1.18e-3 & -7.30e-8 & 1.30e-10 \\
%\end{tabular}
%\caption{Relative change in entropy \eqref{entropy-change} for solutions of the
%Euler equations \eqref{euler} with initial data given by \eqref{smooth-IC}.}
%\label{tab:entropy-change}
%\end{table}

%\begin{table} % Relative entropy change!
%\centering
%\begin{tabular}{l|rr}
%$\Delta \xeul$ & Clawpack & Pseudospectral \\ \hline
%1/16  & 4.79e-3 & 4.71e-8\\
%1/32  & 1.95e-4 & 1.55e-9\\
%1/50  & 2.31e-5 & 1.08e-9\\
%1/100 & 9.23e-7 & 5.94e-10\\
%1/200 & 7.30e-8 & 6.97e-10\\
%\end{tabular}
%\caption{Relative change in entropy \eqref{entropy-change} for solutions of the
%Euler equations \eqref{euler} with initial data given by \eqref{smooth-IC}.}
%\label{tab:entropy-change}
%\end{table}

\begin{figure}
    \centering
    \includegraphics[width=0.7\linewidth]{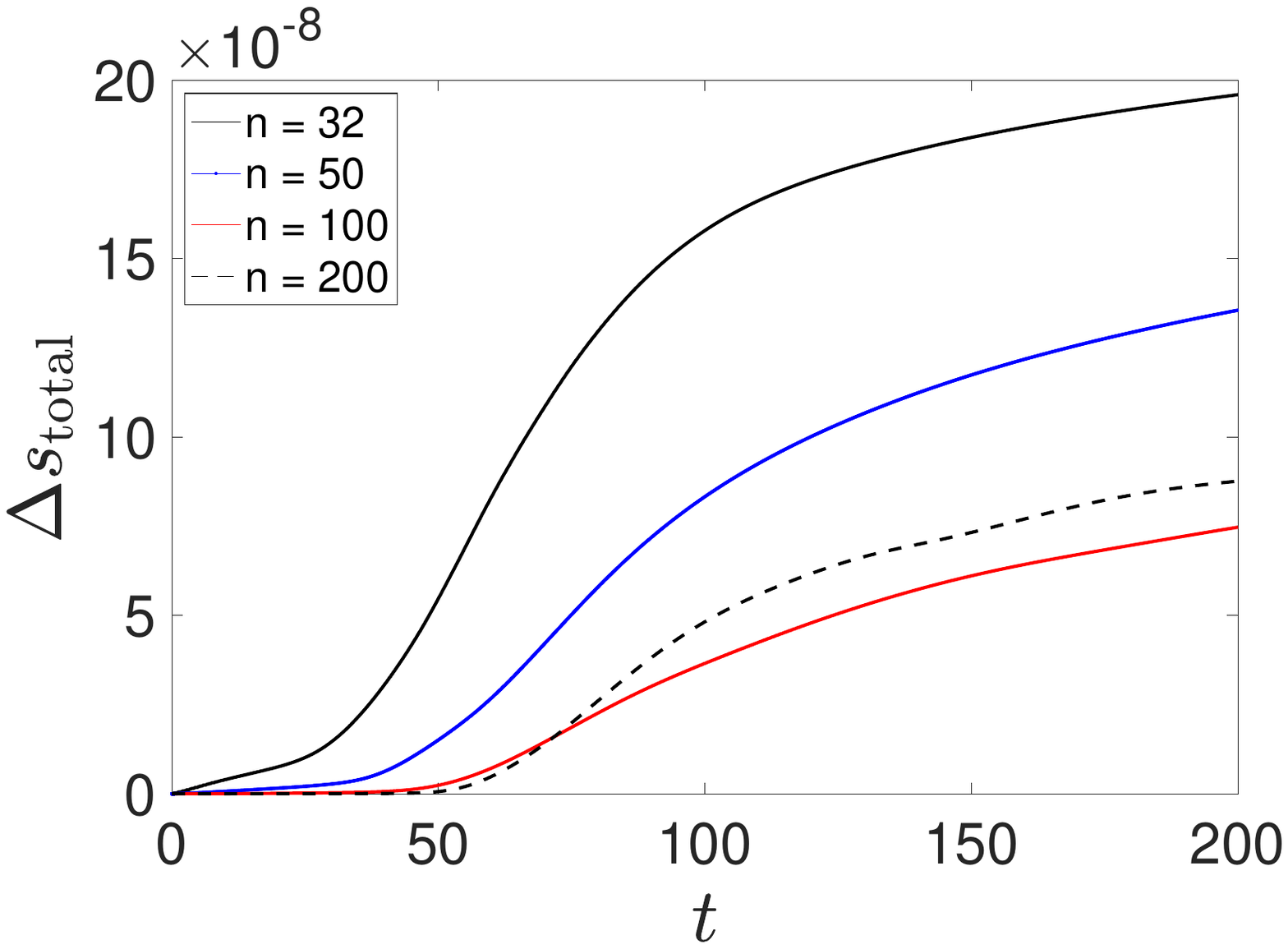}
    \caption{Evolution of the total entropy $\int_0^L s(x,t)\,dx$ as a function of time computed by the pseudospectral method with various number of gridpoints per period.}
    \label{fig:entropy_vs_time}
\end{figure}

\begin{table}
\centering
\begin{tabular}{l|rr}
$\Delta \xeul$ & Clawpack & Pseudospectral \\ \hline
1/16  & 6.02e-1 & 5.93e-6\\
1/32  & 2.45e-2 & 1.96e-7\\
1/50  & 2.90e-3 & 1.36e-7\\
1/100 & 1.16e-4 & 7.47e-8\\
1/200 & 9.18e-6 & 8.76e-8\\
\end{tabular}
\caption{Change in entropy $\Delta S_\textup{total}$ for solutions of the
Euler equations \eqref{euler} with initial data given by \eqref{smooth-IC}.}
\label{tab:entropy-change}
\end{table}
The behavior of the total entropy variation as a function of time is reported in Figure \ref{fig:entropy_vs_time}.
Although we cannot exclude the production of a small amount of entropy, with fully resolved calculations the change in the entropy at final time appears only on the eleventh digit. 

\subsection{A criterion for shock formation}
A careful study of shock formation for the $p$-system in a periodic medium was conducted in
\cite{ketcheson2012shock}, and an empirically-supported criterion for shock formation was
formulated there.  This criterion indicates that in a non-uniform medium a propagating front will form a shock if and only if the speed of the front exceeds a critical threshold, namely, the maximum speed of signal propagation in the state ahead of the front, given by the harmonic average of the sound speed:
\begin{align}
    c_\text{max} & = \avg{c^{-1}}^{-1}
\end{align}
where $c=\sqrt{-p'(v)}$.
Note that the speed of propagation of the front itself is altered by the variation in the background; for details see \cite{ketcheson2012shock}.
%The second quantity is the effective shock speed in the variable
%medium, given by
%\begin{align}
%    w_\text{eff} & = \sqrt{\left\langle\frac{[v]}{[p]}\right\rangle^{-1}},
%\end{align}
%where $[q]=q_r-q_\ell$ is the jump across the shock.
%Note that $w_\text{eff}$ is essentially a harmonic average
%of the Rankine-Hugoniot shock speed.  
%Computational experiments suggest that shock formation occurs when $w_\text{eff}>c_\text{max}$ and is
%avoided otherwise \cite{ketcheson2012shock}.
This criterion cannot be applied precisely to the propagation of localized pulses, since for a pulse the asymptotic left and right states are the same, and the theory was formulated in terms of a propagating front.

Instead, we consider a series of shock-tube-like problems with an oscillating density.  Specifically, we take
\begin{subequations} \label{shocktube-setup}
\begin{align}
    \rho(\xeul,0) & = 1 + \cos(2\pi \xeul) \\
    u(\xeul,0) & = 0 \\
    p(\xeul,0) & \revA{=} \frac{1-\tanh(\xeul/2)}{2} p_\ell + \frac{1+\tanh(\xeul/2)}{2} p_r
    %\begin{cases} p_\ell & x<= 0 \\ 1 & x>0.\end{cases}
\end{align}
\end{subequations}
The pressure is approximately a step function,
taking the value $p_\ell$ for $x<0$ and $p_r$ for $x>0$, but the hyperbolic tangent function is used to smooth it out slightly so that there is not a shock at $t=0$.  In this  way we can assess whether a shock forms later.
We take $p_r=1$ and vary $p_\ell$.
We see that the solution of this problem resembles that of a Riemann problem for a dispersive nonlinear wave equation, with an oscillatory region forming between the left-going rarefaction and right-going front.  Our interest here is in assessing the speed of this front and whether it is truly a shock.  We see clearly that in the first two figures the front is moving slower than $c_\text{max}$, while in the last two figures it is clearly moving faster than $c_\text{max}$.
Thus we expect (based on the theory from \cite{ketcheson2012shock}) shock formation in the last two cases but not in the first two.  The middle cases are more ambiguous.

\begin{figure}
\centering
    \includegraphics[width=\textwidth]{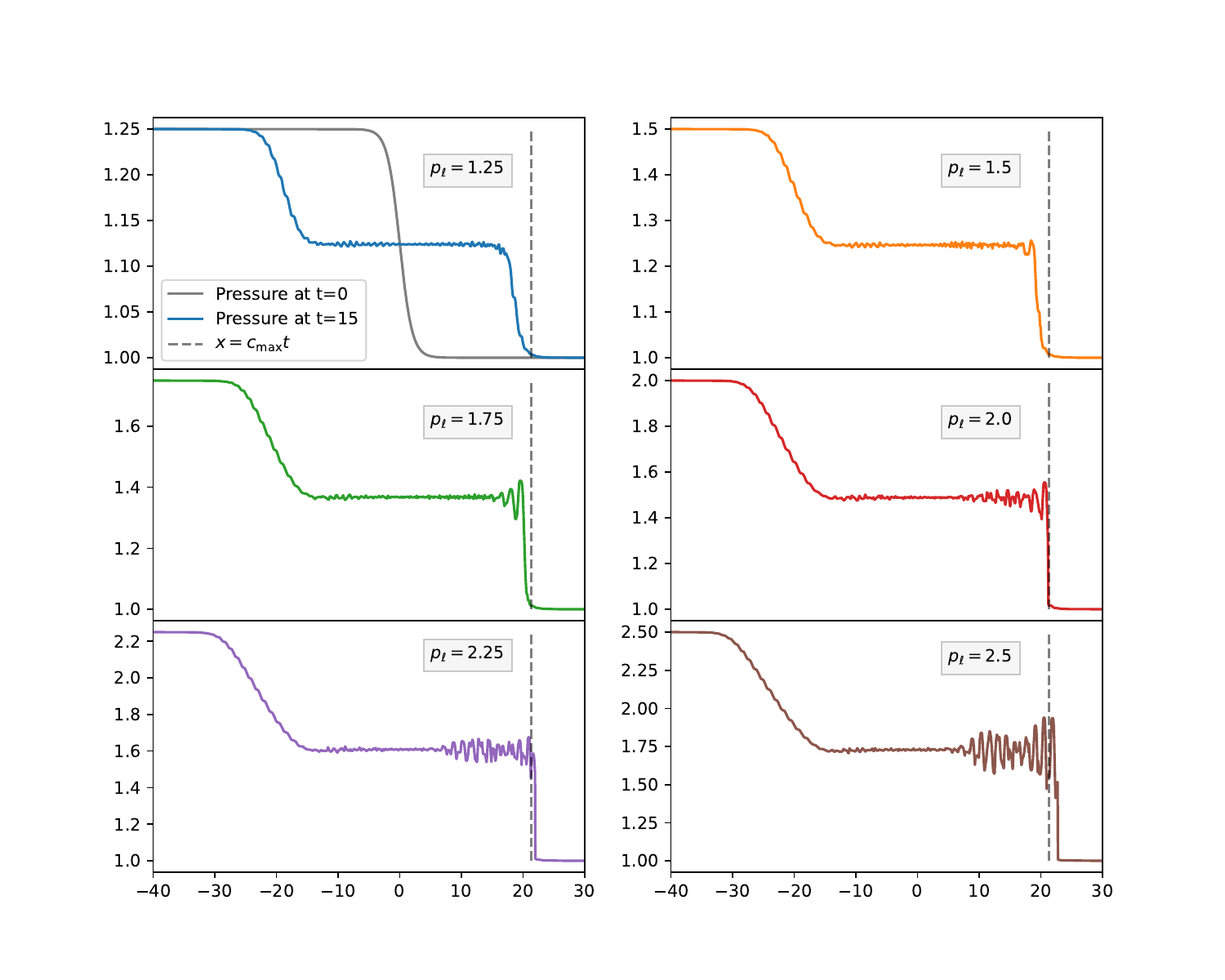}
    \caption{
        Results of variable-entropy shock tube experiments \eqref{shocktube-setup}, with $p_\ell$ ranging from 1.25 to 2.5.  The pressure at $t=15$ is plotted, using colors that correspond to those of Figure \ref{fig:lep}.  The initial pressure is also shown in the top-left plot in grey.  The maximum speed of small-amplitude perturbations in the right state is indicated by the dashed vertical line.
     }
 \label{fig:shocktubes}
\end{figure}

To test this expectation, we compute a measure of the local entropy production (LEP):
\begin{align}
    \lep^n_j & =\frac{S^{n}_j-S^{n-1}_j}{\Delta t} + \frac{\psi^n_{j+1}-\psi^n_{j-1}}{2\Delta \xeul}.
\end{align}
Local entropy production a very useful tool to detect the presence and nature of possible singularities in the numerical solution of quasilinear hyperbolic systems of conservation laws. 
This tool was originally introduced by G. Puppo 
\cite{puppo2002numerical,puppo2004numerical} in the context of finite volume central schemes, and later extended to other schemes. 
In \cite{puppo2011numerical} the authors studied the dependence of this quantity on the spatial mesh size $h$. When applied to the Riemann problem for the Euler equations, they proved that when using a method of order $p$ in space and time, the local numerical entropy production 
scales as $\order(h^p)$ in the smooth regions, and as $\order(h^{-1})$ near the shock. They also observed that LEP scales as $\order(h)$ at the rarefaction corner, as $\order(h^0)$ at the contact. In the same paper
the LEP has been adopted
as an indicator for adaptive mesh refinement in finite volume methods on uniform grids.
In in \cite{puppo2016well} and  \cite{semplice2016adaptive} 
indicators based on LEP have also been adopted in the case of non uniform grids.

\begin{figure}
\centering
    \includegraphics[width=0.6\textwidth]{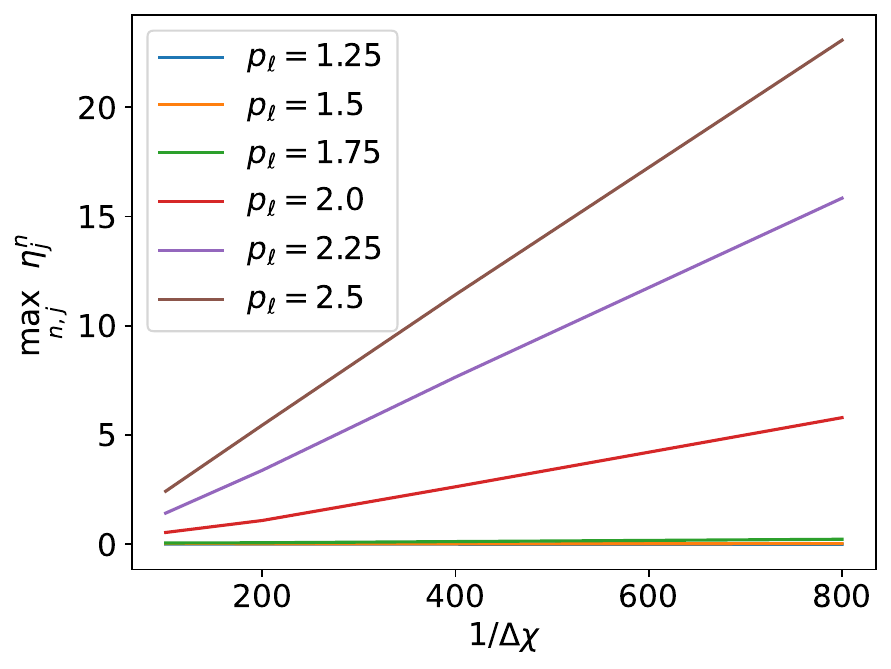}
    \caption{
        Maximum local entropy production versus mesh spacing, for shock tube simulations with varying values of $p_\ell$.
        Line colors correspond to those of Figure \ref{fig:shocktubes}.
        For larger values, the entropy production grows linearly as the mesh size is decreased, indicating the presence of one or more shocks.  For smaller values of $p_\ell$
        this measure suggests that shocks are not present or are very weak (note that the plots
        for $p_\ell=1.25, 1.5, 1.75$ lie almost on top
        of one another).
     }
 \label{fig:lep}
\end{figure}

In Figure \ref{fig:lep} we plot the maximum (over time and space) of the absolute value of $\lep^n_j$ for each shock tube experiment, for various values of $\Delta \chi$ and $p_\ell$.
For values $p_\ell\ge 2$, this value
grows linearly as the mesh is refined, indicating the presence of one or more shocks.  For smaller values of $p_\ell$ this measurement suggests that shocks are not present or are extremely weak.

%However, we can apply it in an approximate manner by
%taking the peak of the initial pulse as $q_\ell$ and the ambient state as $q_r$.  Doing so, we find that
%the results above agree roughly with this theory:
%we have $w_\text{eff}<c_\text{max}$ for the scenario studied in Figure \ref{fig:compare3} but $w_\text{eff}>c_\text{max}$ for the scenario in Figure \ref{fig:compare3_shock}.
%For a more thorough and detailed study of shock formation in the $p$-system
%with periodic coefficients we refer the reader to \cite{ketcheson2012shock}.

\section{Perturbations in a random quasi-periodic background}\label{sec:random}

%\red{In this section we try to give a plausibility argument 
%that explains how acoustic waves propagating into a constant background should break into a shock, still they do not break for non uniform, non periodic background. 
%We apply the ideas above to a very simple model 
%for acoustic waves in the atmosphere.  In the 1D Euler model,
%any acoustic perturbation will eventually break, but small perturbations break only after a long time.  We work out the breaking time in Section \ref{sec:acoustic-breaking}, and perform a simulation of waves for Euler equations on an almost periodic background in the next section.}

Here we have observed that genuinely nonlinear waves do not break up into shocks, and presumably maintain the regularity of the initial conditions, provided they propagate on a periodic background. 
Similar observations have been made for other hyperbolic systems; see e.g. \cite{KLR2023,busaleh2024homogenized}.
A natural question arises: how fundamental is the periodicity of the background? What if the background is oscillatory but not strictly periodic? 

\begin{figure}
    \centering
    \includegraphics[width=\textwidth]{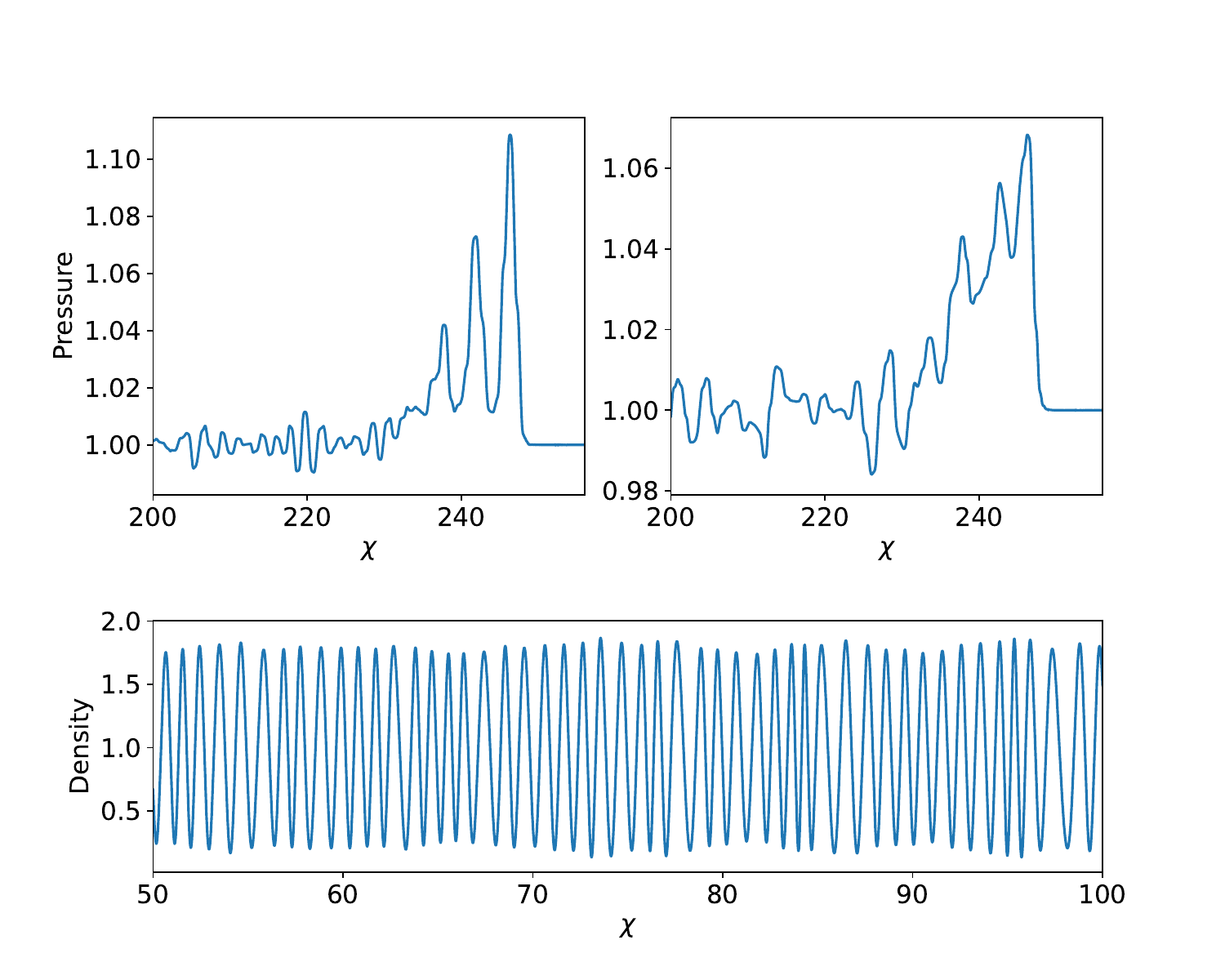}
    \caption{Propagation of a pressure perturbation on a random background.  The upper figures show results of two different realizations, one obtained by using the algorithm described in Appendix \ref{sec:random-background} using 320,000 smoothing iterations (left panel) and the other obtained by the same algorithm, using 20,000 iterations. 
     The lower figure shows an example of the randomly oscillating density field.  For visual clarity, only a subset of the domain is shown. 
 The relative change in entropy for the two simulations shown is $-2.37\times10^{-7}$ and $-1.78\times10^{-7}$, respectively.}\label{fig:random}
\end{figure}

Here we conduct some very simple and limited experiments to begin probing the answer to this question.
We consider a background density of the form 
\[
    \rho(\xeul,t=0) = 1 + 0.8 A(\xeul)\sin(2\pi B(\xeul) \, \xeul)
\]
where $A(\xeul)$ and $B(\xeul)$ are random functions centered about unity (see Appendix \ref{sec:random-background}).  We conduct two realizations of these random functions, with differing amounts of variation (specifically, the number of smoothing iterations is set to 320,000 for the smoother one and 20,000 for the less-smooth one).  We apply the high-order SharpClaw FV algorithm from Clawpack, solving on the domain $[-256,256]$ up to $t=200$, just as in the previous section.  We use a grid with $\Delta \xeul = 1/200$.  Solutions are shown in Figure \ref{fig:random}, where we have zoomed in to show the leading part of the wave.
A part of one of the random entropy profiles is also shown.  We see that the solution no longer consists of such neat solitary waves, although a roughly similar pattern emerges (with a lot more small-scale oscillations).

\begin{figure}
    \centering
    \includegraphics[width=0.5\textwidth]{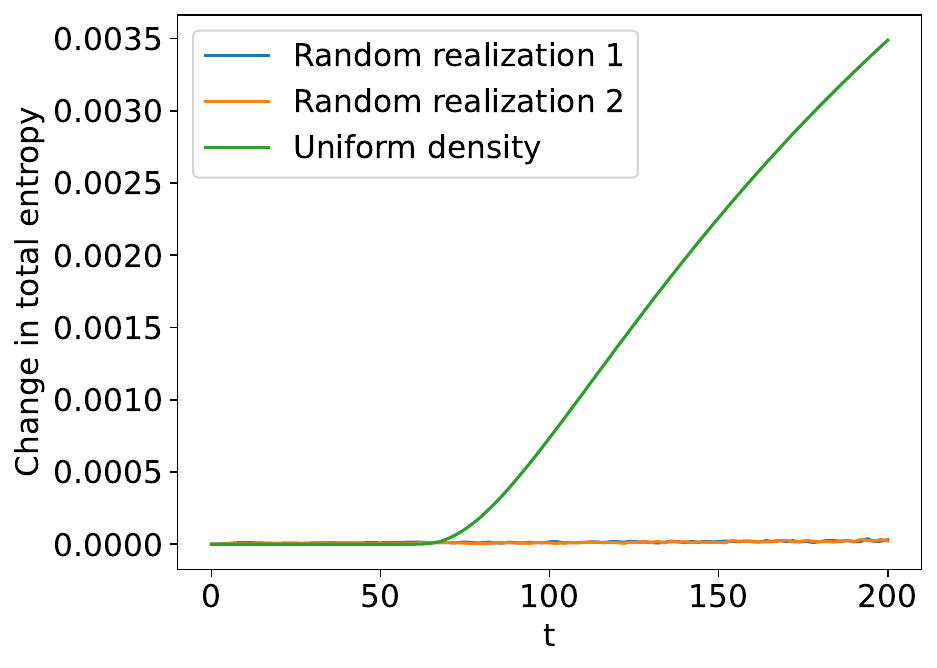}
    \caption{Entropy change over time for two realizations of a random oscillating density field, compared with a spatially uniform density field.}\label{fig:ent_over_time_random}
\end{figure}

In Figure \ref{fig:ent_over_time_random} we show the change in
total entropy over time for the solutions over a random field compared to the same initial pressure perturbation but with a constant (in space) density field.  It is clear that the entropy change in the presence of a random oscillating density field is orders of magnitude smaller, and if shocks form then they are extremely weak.  We also computed the relative change in entropy (eq. \eqref{entropy-change}) for the two random realizations; for the smoother density field we obtain $1.24\times 10^{-7}$ and for the less-smooth density field we obtain $-1.78 \times 10^{-7}$.  These values are of the same order of magnitude (though somewhat larger) as what was obtained for a strictly periodic density field with the same grid spacing (see Table \ref{tab:entropy-change}).  This further supports the hypothesis that little or no shock formation occurs in these scenarios.
Similar results (not shown here) were obtained with additional realizations of the random density field.

This experiment suggests that periodicity of the background is not a necessary condition in order to have a smooth solution for long time. 
%The result is in agreement with the analysis performed by Temple and Young in \cite{temple2023nonlinear}. 
Further investigation in this direction, in order to understand
the conditions on the background under which genuinely nonlinear waves of quasilinear systems remain smooth for arbitrarily long times represents an interesting and challenging problem for future study.

%\begin{figure}
%\centering
%    \includegraphics[width=0.48\textwidth]{Figures/pressure_zoom.pdf}
%    \includegraphics[width=0.48\textwidth]{Figures/pressure123.pdf}
%    \includegraphics[width=0.48\textwidth]{Figures/pressure4569.pdf}
%    \includegraphics[width=0.48\textwidth]{Figures/pressure74391.pdf}
%    \includegraphics[width=0.96\textwidth]{Figures/density1.pdf}
%    \caption{Upper panels: signal propagation on a perfectly periodic background (upper left panel) and on a weakly modulated background (upper right and middle panels). 
%    Three different seeds have been adopted to obtain three different random sequences. If the perturbation is small the qualitative behaviour is similar. Both numerical solutions are obtained by the spectral method with Runge-Kutta 4, using a CFL number 0.9. Lower panel: the weakly modulated (almost periodic) initial density, corresponding to the upper right panel, in an certain interval. }
%    \label{fig:non_periodic}
%\end{figure}

\section{Conclusion}
We have shown that a pressure perturbation propagating on a background of a gas with
periodically (or almost periodically) varying density (or equivalently, entropy) undergoes an effective dispersion
that can lead to the formation of solitary waves that seem to persist for arbitarily long times.  This behavior can be accurately
described by a pair of high-order constant-coefficient PDEs obtained via perturbation theory.
\revA{The present work seems to be the first time that this behavior has been observed for the compressible Euler equations on a non-periodic domain.}
In agreement with  recent work on nonlinear acoustics \cite{temple2023nonlinear}, our
results suggest that in general solutions of the compressible Euler equations need not
decay to a constant state.  We note that whereas the results in \cite{temple2023nonlinear} deal only with very small amplitude (nonlinear) perturbations, our results suggest the existence of
non-decaying solutions of relatively large amplitude.

The behavior observed here seems to be intimately connected to work on resonant wave interactions \cite{majda1984resonantly,majda1988canonical,pego1988some,hunter2019resonant}.  Although we have not explored this aspect here, the traveling wave solutions we observe are composed of variation in all three characteristic fields (cf. \cite{leveque2003phase}).  It is remarkable that such interactions persist indefinitely on an unbounded domain.

We also observe that, given the regularity of the solution, one can effectively adopt simpler and more accurate methods for the numerical solution of the Euler system, such as, for example, a pseudospectral method coupled with a standard ODE integrator, with great gain in computational efficiency.  Of course, this requires advance knowledge of the nature of the solution, for example that is is sufficiently regular to be treated by pseudsospectral methods, 
which can be facilitated partially by the shock formation criterion proposed in \cite{ketcheson2012shock} and revisited here in Section \ref{sec:entropy}.  Nevertheless, there is 
not yet a complete theory to guarantee the avoidance of shock formation for general classes of initial data.

\revA{We have not attempted herein to provide a quantitative estimate
of the error in our approximation or of the regime in which it is valid, although the numerical solutions we have studied provide some limited data regarding this question.  An answer to this question for some linear wave equations can be found in \cite{allaire2022crime}.
For nonlinear equations like those studied herein, this remains an
interesting and challenging area for future work.
}

Many interesting open questions are raised by this work.  For instance, is it possible
to prove that there exist large-amplitude non-breaking solutions of the 1D Euler equations,
and is there a limit to how large they can be?  %What are the properties of these 
%quasi-traveling wave solutions, in terms of their shape, speed, and interactions?
Could these waves be generated and observed experimentally?

At this point it is natural to expect that similar behavior to what is observed here
(and recently in the shallow water system \cite{KLR2023}) may arise in other
hyperbolic systems with spatially-periodic structure.  Based on the work of
Temple \& Young \cite{temple2009paradigm,temple2011time,temple2023nonlinear}, this behavior may require the presence of two nonlinear
characteristic fields and one linearly degenerate field.  An investigation of
the necessary and sufficient hyperbolic structure for such behavior to arise is
the subject of ongoing work.

\appendix
\section{Appendix}\label{sec:appendix}

\subsection{Linear stability analysis of the solution of system \eqref{homog_delta4}} \label{A:stability}
Here we prove that system \eqref{homog_delta4} is always linearly unstable, for all $k$, and that the initial value problem for the linearized system is ill-posed. 

Let us consider the linearized system \eqref{homog_delta4}, and neglect the terms on the right hand side. 
The system can be written as 
\begin{subequations} \label{lin_ttt}
\begin{align} 
    p_t + \beta_1 u_x -\beta_2 \partial^3_t p + \beta_3 \partial^5_t p & = 0\\
    u_t + p_x = 0
\end{align}
\end{subequations} 
where we set
\[
    \beta_1 = \frac{G(p_*)}{\mean{K^{-1}}}, \quad 
    \beta_2 = \delta^3\mu\frac{\mean{K^{-1}}}{G(p_*)}, \quad
    \beta_3 = \delta^4\frac{\zeta}{\mean{K^{-1}}} \alpha_7.
\]
We assume $\beta_1>0$, $\beta_2 > 0$, while $\beta_3\geq 0$ ($\beta_3 = 0$ corresponds to the model which neglects terms of order higher than $\order(\delta^2)$).
We look for solution of the form 
\[
    p = \hat{p}\,e^{i(kx-\omega t)},\quad
    u = \hat{u}\,e^{i(kx-\omega t)}.
\]
After dividing by the imaginary unit $i$ and collecting the terms in $\hat{p}$ and $\hat{u}$ we obtain the homogeneous system
\begin{align*}
    (-\omega - \beta_2\omega^3-\beta_3\omega^5) \, \hat p   +  \beta_1 k \,\hat u & = 0,\\
    k \, \hat p   -  \phantom{\beta_,} \omega \,\hat u & = 0 .
\end{align*}
Non trivial solutions of the system exist if the determinant of the coefficient matrix is zero. This gives the dispersion relation
\[
    \omega^2 + \beta_2\omega^4 + \beta_3 \omega^6 - \beta_1 k^2 = 0.
\]
Let $Y=\omega^2$. Then $Y$ satisfies the cubic equation 
\begin{equation}\label{A:cubic}
    f(Y) = \beta_1k^2, \quad f(Y)\equiv Y + \beta_2 Y^2 + \beta_3 Y^3.
\end{equation}
For $\beta_3 = 0$, $Y$ satisfies the equation
\[
   \beta_2Y^2 + Y -\beta_1k^2 = 0
\]
which has two real roots of opposite sign, say $Y_+>0$, $Y_-<0$. 
The equation $\omega^2 = Y_-$ has two imaginary roots, one of which with positive imaginary part, corresponding to an unstable mode which grows exponentially. 

If $\beta_3>0$, since the function $f(Y)$ is monotonically increasing and $f(0)=0$, there will be only one positive root, which corresponds to two real values of $\omega$. The other roots of \eqref{A:cubic}
are either negative real or complex. Let us denote by $Y_*$ one of such roots. 
If $Y_*<0$ then $\omega=i\sqrt{-Y_*}$ corresponds to an unstable mode. If $Y_*$ is a complex root, one of the two roots of $\omega^2 = Y_*$ will have positive imaginary part, therefore corresponding to an unstable mode. 

Finally, as $|k|\to \infty$, so the corresponding roots diverge, i.e.\ $|\omega|\to\infty$, making the initial value problem ill-posed.

\subsection{Linear stability of LeVeque \& Yong system}
In \cite{leveque2003}, the homogenized equations obtained there
are converted to a form in which only x-derivatives appear
(except for the evolution terms).  
Here we derive the linear stability relation for
the analogous system, taking the Euler equation of state rather than the exponential stress-strain relation used there.
Using \eqref{correspondence}
we obtain the system
\begin{align}
    p_t + \frac{G}{\mean{K^{-1}}} u_x & = \frac{\delta^2}{\mean{K^{-1}}} \left(-\mean{K^{-1}} \mu G u_{xxx} + 2 \mu G' p_x u_{xx} - 
    \mu G'' (p_x)^2 u_x \right) + \order(\delta^3) \\
    u_t + p_x & = 0.
\end{align}
The linearized system around state $(u_*=0, p_*)$ reads
\begin{align}
    p_t + \beta_1 u_x +\tilde{\beta}u_{xxx} & = 0 \\
    u_t + p_x & = 0
\end{align}
where $\beta_1$ is defined in the previous section, and $\tilde{\beta} = \delta^2\mu G(p_*) > 0$.
We look for solutions of the form 
\[
    p = \hat{p}\,e^{i(kx-\omega t)},\quad
    u = \hat{u}\,e^{i(kx-\omega t)}.
\]
After dividing by the imaginary unit $i$ and collecting the terms in $\hat{p}$ and $\hat{u}$ we obtain the homogeneous system
\begin{align*}
    (\beta_1 k-\tilde{\beta} k^3) \hat{u}-\omega \hat p & = 0\\
    -\omega \hat u + k \hat p & = 0
\end{align*}
Non trivial solutions of the system exist if the determinant of the coefficient matrix is zero. This gives the dispersion relation
\[
    \omega^2 =   \beta_1 k^2 -\tilde{\beta} k^4,
\]
which gives the two branches
\[
    \omega = \pm |k|\sqrt{\beta_1-\tilde{\beta}k^2}
\]
therefore the linear modes are unstable for $|k|>\sqrt{\beta_1/\tilde\beta}$, leading to ill-posed initial value problems.

\subsection{Randomly modulated background}
\label{sec:random-background}
Here we describe how we construct the coefficients $A(x)$ and $B(x)$ which modulate the amplitude and space frequency of the originally sinusoidal background. We describe the construction of $A(x)$. 
First we solve the stochastic differential equation 
\[
    dA(x) = - K_A A(x) \, dx + \sigma_A (\xi-0.5) \, \sqrt{dx}
\]
numerically on the grid which discretizes the domain $[-L,L]$, with a grid f stepsize $dx = 2L/N$, starting from $A(-L) = 0$. Here $K_A$ and $\sigma_A$ are two non-negative constants, and $\xi$ denotes a random number with uniform distribution in $[0,1)$. In this
way we obtain a discretization of a H\"older continuous function with exponent $1/2$, with zero expected value. 
Then we smooth the obtained discrete values $A_i\approx A(x_i),i=1,\ldots,N$ by discrete diffusion equation, i.e. 
we iterate
\[
    A_i^{(n+1)} = A_i^{(n)} + 0.25(A_{i+1}^{(n)}-2A_i^{(n)}+A_{i-1}^{(n)}),\> i=1,\ldots,N,
\]
imposing periodic boundary conditions, and with $n=1,\ldots,N_{\rm iterations}$.
At the end of the process we add 1 to $A$, so that its expected value is 1. Similarly we act for $B(x)$.
Specifically, in our simulation we adopted $K_A = K_B =1$, $\sigma_A = 0.24$, $\sigma_B = 0.015$.
In order to make the result reproducible we set a seed of the random number generator with the Matlab command {\tt rng}.

\section*{Acknowledgments}
The authors would like to thank Giuseppe Virgilio Minissale for pointing out the approximation that allows to use Petviashvili technique to compute the traveling wave. 
G.~Russo would like to thank the Italian Ministry of University and Research (MUR) for the support of this research with funds coming from PRIN Project 2022 (N.\ 2022KA3JBA entitled  ``Advanced numerical methods for time dependent parametric partial differential equations with applications'') and King Abdullah University of Science and Technology (KAUST) for hosting him during part of the time this work was completed.  D. Ketcheson thanks KAUST for the funding that supported this work.
\bibliographystyle{plain}
\bibliography{references}

\end{document}